\documentclass[onefignum,onetabnum]{siamart171218}
\usepackage[framed,numbered,autolinebreaks,useliterate]{mcode}
\usepackage{lipsum}
\usepackage{amsfonts}
\usepackage{graphicx}
\usepackage{epstopdf}
\usepackage{algorithmic}
\ifpdf
  \DeclareGraphicsExtensions{.eps,.pdf,.png,.jpg}
\else
  \DeclareGraphicsExtensions{.eps}
\fi

\newsiamremark{remark}{Remark}
\newsiamremark{hypothesis}{Hypothesis}
\crefname{hypothesis}{Hypothesis}{Hypotheses}
\newsiamthm{claim}{Claim}

\headers{A Domain Overlapping Algorithm with Nonlinear Mapping for Collocation-Based Solutions of  Eigenvalue Problems}
{ Jinwei Yang and Vinod Srinivasan}

\title{A Domain Overlapping Algorithm with Nonlinear Mapping for Collocation-Based Solutions of  Eigenvalue Problems\thanks{Submitted to the editors DATE.
\funding{This work was funded by NSF under the number}}}

\author{Jinwei Yang\thanks{Department of Mechanical Engineering, University of Minnesota, Minneapolis, Minnesota, USA}
\and Vinod Srinivasan \thanks{Department of Mechanical Engineering, University of Minnesota, Minneapolis, Minnesota, USA (\email{vinods@umn.edu}).}}

\usepackage{amsopn}

\makeatletter
\newcommand*{\addFileDependency}[1]{
  \typeout{(#1)}
  \@addtofilelist{#1}
  \IfFileExists{#1}{}{\typeout{No file #1.}}
}

\makeatother
\newcommand*{\myexternaldocument}[1]{%
    \externaldocument{#1}%
    \addFileDependency{#1.tex}%
    \addFileDependency{#1.aux}%
}

\ifpdf
\hypersetup{
  pdftitle={A Domain Overlapping Algorithm with Nonlinear Mapping for Collocation-Based Solutions of  Eigenvalue Problems}, pdfauthor={Jinwei Yang and Vinod Srinivasan}}\fi

\myexternaldocument{ex_supplement}
\usepackage{ulem}
\usepackage{soul}
\begin{document}

\maketitle
\begin{abstract}
In this paper, we present a set of domain decomposition algorithms incorporating nonlinear mapping to enhance the numerical solution of eigenvalue problems, particularly for sharp-interface geometries. Traditional spectral collocation methods using Chebyshev polynomials exhibit challenges in handling discontinuities and steep gradients within the computational domain. To address these issues, we introduce four domain overlapping algorithms that improve spectral convergence by strategically distributing grid points while maintaining the smoothness of higher-order derivatives. The proposed methods are demonstrated on the 1D Burgers equation and the eigenvalue problem for instability in laminar 3D channel flow with viscosity stratification. Comparisons with conventional numerical techniques, including global mapping and Chebfun’s splitting method, illustrate the advantages of our approach in terms of accuracy and convergence. The effectiveness of the proposed methods in capturing complex physical phenomena with reduced computational cost makes them promising tools for solving a broad class of problems involving sharp interfaces and multi-domain structures.
\end{abstract}



\section{Introduction}
In many problems defined on finite domains, such as wall-bounded fluid flows, the steepest gradients in the quantity of interest occur near the boundaries of the domain. Chebyshev spectral collocation techniques are attractive for such situations, because of their high accuracy, exponential convergence properties, as well as a grid point distribution that clusters points near the regions of the sharpest gradients \textbf{ref}. When the problem is defined on other geometries, such as semi-infinite or infinite domains, this technique can be used with adaptation, such as by using a set of bases functions that satisfy the required exponential decay properties \cite{boyd1987orthogonal, boyd1987spectral}, or by using suitable mappings that transform the physical space into the Chebyshev domain. However, when sharp or weak discontinuities exist in the interior of a domain, then it is not easy to find mappings that will result in an adequate point clustering near the region of discontinuity. Implementing a collocation method with such a mapping leads to loss of the smooth convergence properties of Chebyshev polynomials, and/or requires a very high number of grid points in the entire domain when appropriate mappings are found, it may still lead to complexity of the code.

To avoid these issues, methods involving domain decomposition have been proposed. Two families of multi-domain methods are in use today: One is based on the decomposition of the domain into overlapping subdomains, with information being transferred from one subdomain to another through the overlapping regions. All of those techniques come from the Schwarz method which was first proposed to prove the existence of the solutions for elliptic boundary value problems. More examples of these algorithms are given in the literature [5, 9]. The second family of methods involves decomposing the physical domain into adjoining subdomains by using suitable conditions at the interfaces \textbf{ref}. 

In recent years, the MATLAB suite of codes called Chebfun has attracted significant interest, due to its ability to provide very accurate automatic solutions to differential equations\cite{driscoll2014chebfun}. Chebfun, an open-source project, also includes a domain decomposation function, called the split method, which creates piece-wise polynomial approximations\cite{pachon2010piecewise}. When splitting is enabled, if a Chebyshev interpolant is unable to represent the function accurately at a specified maximum degree on an interval, the interval is bisected; this process is recursively repeated on the subintervals. A drawback of Chebfun's splitting approach is that the resulting representation does not ensure anything more than $C^0$ continuity. In other words, for a function that is analytic on and around an interval, Chebyshev polynomial interpolation provides spectral convergence. However, if the function has a singularity close to the interval, the rate of convergence is near one \cite{aiton2018adaptive,aiton2019adaptive,aiton2020preconditioned}. A method to remove this issue is to split the domain into overlapping intervals and use an infinitely smooth partition of unity to blend the local Chebyshev interpolants \cite{aiton2018adaptive}. However, this method can not cluster the node distribution as required for some situations where a sharp discontinuity exists, or when steep gradients ar present. In the present case, we propose a new overlapping domain decomposition algorithm which ensures smoothness of higher-order derivatives near the region of discontinuity. The premise of the present study is that by distributing points in a specific manner on both sides of an interface of zero or finite but small thickness, one can progressively improve the rate of convergence. We first show that a one-point overlap (i.e. standard domain decomposition)) can partially resolve the issue, by using information from the two adjoining subdomains. When compared to the case of a single global mapping over the entire domain, however, it can lead to large errors at the overlap point in some cases due to poor control over point distribution.  Next, by using a two-point overlapping technique, where the points from the two domains overlap with each other, one can show significantly reduced error and faster convergence. Finally, we generalize this approach and explore the effects of adding more overlapping points in a region of steep gradients, while relaxing the requirement of exactly overlapping grid points.

Coordinate transformations are very important in spectral methods since the Chebyshev polynomial to cosine change of variable greatly simplifies computer programs for solving differential equations. More details about nonlinear mapping are available in   \cite{boyd2001chebyshev}. Examples of infinite and semi-infinite intervals mapped to $[-1,1]$ will be illustrated in this paper, followed by the use of Chebyshev polynomials. Generally, when the computational interval is unbounded, a variety of options are available, and only some of them require a change of coordinate. For example, on the infinite interval [$-\infty$, $\infty$], Hermite functions or sinc functions are good basis sets. On the semi-infinite domain, [0, $\infty$], the Laguerre functions give exponential convergence for functions that decay exponentially to infinity. A second option is domain truncation, in which the problem is solved on a large but finite interval, $[-L,L]]$ (infinite domain) or $[0,L]$ (semi-infinite) using Chebyshev polynomials with appropriate;y mapped argument respectively. The third option, and the one we will use in this paper is to use a nonlinear mapping that will transform the unbounded interval into [-1,1], so that we can apply Chebyshev polynomials without truncating the computational interval to a finite size. An early paper by Grosch and Orszag \cite{grosch1977numerical} compared the logarithmic map and algebraic map and they found that both maps work well for functions that decay exponentially in space. Boyd \cite{boyd1987spectral, boyd2001chebyshev} further introduced two families of transformed Chebyshev functions for the infinite and semi-infinite domains, that satisfy natural boundary conditions at infinity. 

The paper is organized as follows. The proposed overlapping grid algorithm is introduced in \cref{sec:alg}. Code validation and numerical results for various values of number of overlapping points are given in \cref{sec:experiments}, and the conclusions follow in \cref{sec:conclusions}.

\section{Proposed domain overlapping algorithm}
\label{sec:alg}

We use Chebyshev interpolants for our method because they enjoy spectral convergence. Further, as explained above, they can easily be used with appropriate transforms for unbounded intervals. 
\subsection{Algorithm 1: One point overlapping method}

Let us suppose that the calculation domain is $[0, L]$, $L>1$, with a discontinuity located at $y=L_1$.
One conventional approach that has been previously adopted would use a single nonlinear mapping to cluster points in the vicinity of the discontinuity. For sufficiently sharp gradients in the function, this `global' mapping approach would require a large number of points $N$ to ensure that errors remain under a specific level, since clustering points near the gradient region also increases the overall number of points elsewhere unless an extremely appropriate mapping can be found. The alternative approach that has been more recently followed, as exemplified by chebfun with splitting enabled, is to employ domain decomposition/substructuring.  This allows for a block diagonal formulation of the problem in which fewer points are required for subdomains with smooth variation, and higher order representations (more points) for regions with steep gradients. The piecewise representations of the unknown function in the various subdomains are connected through explicit patching of the function value, sometimes along with the first derivative. Thus, the resulting solution has discontinuous higher-order derivatives, which often may not be problematic.  

The basic premise of the one-point overlap technique is to employ domain decomposition, in which only the boundary point is common, but to consider the boundary point as being part of both subdomains, instead of being separate as is done in methods with patching. For a domain with a single interior discontinuity, we split the global domain into two subdomains with a common boundary at $y=L_1$.  using appropriate mappings that cluster points near the boundaries of these subdomains, i.e. in the vicinity of the discontinuity. Enforcing the governing equations at each grid point in the two subdomains using the Chebyshev expansion at the common boundary (the `overlap' point) by using polynomial representations from both subdomains yields an extra equation that needs to be reconciled, in order to produce a non-singular matrix demanded by the solver. 

We illustrate the approach using a linear boundary value problem defined on $y \in [0,\infty]$. 
$$
U^{\prime \prime}(y)=\frac{1}{\theta^2}\left[U(y)-U(y)^2\right](2 U(y)-1)
$$
where the boundary condition is $y(0) = 1$ and $y(+\infty) = 0$.  The equation has a parameter $\theta$ which controls the width of the region over which $U(y)$ changes its value from 1 to 0. The analytical solution 
is $y=\frac{1+\tanh \frac{y-1}{2\theta}}{2}$ and its first derivative at y=1 is determined by $\theta$.

The matrix formulation is outlined below. Suppose $D_1$ and $D_2$ denote the first and second-order derivative matrix operators for the standard Chebyshev method on domain $x \in [-1,1]$,with  the velocity derivative denoted as $\frac{dU}{dy}=D*U$ numerically.  We first truncate the domain $[0,\infty]$ to $[0, L]$ for sufficiently large L and further split this into two parts $[0,L_1]$ and $[L_1,L]$, respectively. For each subdomain, we can introduce nonlinear mappings as follows. 
For the first domain $y_1 = [0,L_1]$, we use 
\begin{equation}
y_1 = \frac{R_c L_1(1-x)}{2 R_c+L_1(1-x^2)},
\end{equation}
which maps the physical domain  $[1,-1]$ to $[0,L_1]$. Here $R_c$ is a free parameter that can be used to adjust the node distribution on this domain. More generally, one can use any reasonable mapping with desired smoothness,  and there are no constraints on the choice of nonlinear mapping. 
For the second domain $[L_1, L]$, we use \textbf{ a linear mapping as follows, though a nonlinear mapping could be used as well}:
\begin{equation}
y_2 = \frac{L_1-L}{2}x + \frac{L_1+L}{2},
\end{equation}
The first-order derivative matrix for the two domains can be derived via the chain rule as follows\cite{boyd2001chebyshev}
\begin{equation}
d_1 = -\frac{[2R_c+L_1(1-x^2)]^2}{R_c L_1 (2R_c+R_1(x-1)^2)}D,
\end{equation}

\begin{equation}
d_2 = \frac{2}{L_1+L_2}D,
\end{equation}

\begin{figure}[h]
\centering
\includegraphics[height=2.0in]{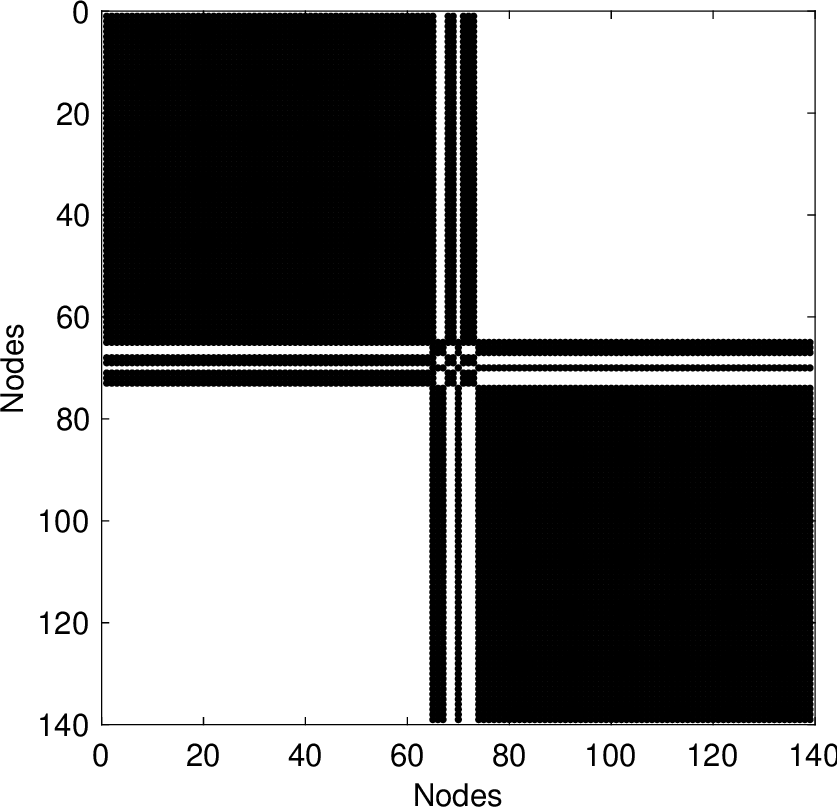}
\includegraphics[height=2.0in]{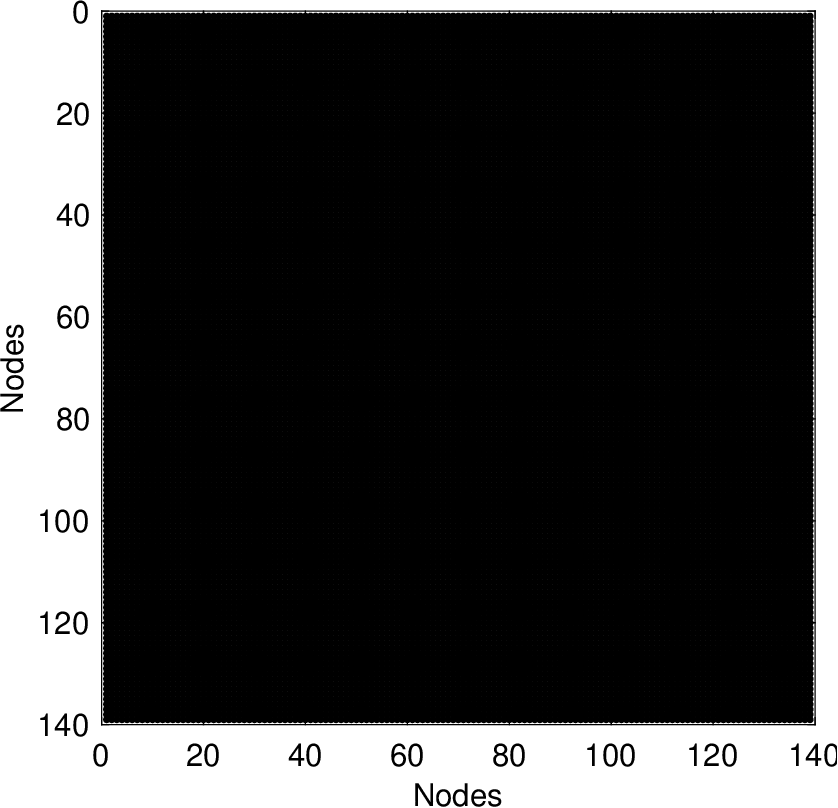}
\caption{Non-zero matrix elements distribution of the semi-global method: (1) The first order derivative matrix and (b) the second order derivative matrix.}
\label{fig:d1d2}
\end{figure}

By using the matrices for each subdomain, we can assemble the global matrix $\hat{D}$. If $u_1$ and $u_2$ are the function profiles in each domain, then one can write:
\begin{equation}
\left[\begin{array}{l}
\frac{d u_{1}}{d y_{1}} \\
\frac{d u_{2}}{d y_{2}}
\end{array}\right]=\left[\begin{array}{lll}
d_{1} & 0\\
0 & d_{2}
\end{array}\right]\left[\begin{array}{l}
u_{1}(y) \\
u_{2}(y)
\end{array}\right] 
\end{equation}
Note that the last row of $d_1$ and the first row of $d_2$ are being written for the common boundary point at $y=L_1$. The first and second-order derivative matrices are shown schematically in Fig. \ref{fig:d1d2}.

 It can be observed from the figure that enforcing the governing equation for the overlapping point leads to an extra equation at 
 $y=L_1$ due to the function at this location having two representations, one from each subdomain.  This needs to be reconciled to obtain a square global matrix. We use the average values of the last row of matrix $d_1$ and the first row of matrix $d_2$ to replace the row corresponding to $y=L_1$. The second order derivative is given by $\hat{D_2}$ = $\hat{D}^2$. Instead of using equal values of 0.5 for the weights of these two rows. In principle, one could use other weights for these rows that sum to unity, though results are known to deteriorate when straying too far from equal weighting [ref]. Examining the matrix in Fig. \ref{fig:d1d2}, it is apparent that this introdices a change in matrix structure from the block diagonal format, adding a single row and column of non-zero elements that use information from all grid points in the domain, and not just from within each subdomain. This is in contrast to \textit{chebfun} with splitting enabled, which appears to use a piecewise continuous representation of the function without passing information from the entire domain to each subdomain. We therefore refer to the proposed algorithm as a `semi-global' method.

\begin{figure}
\centering
\includegraphics[height=2.0in]{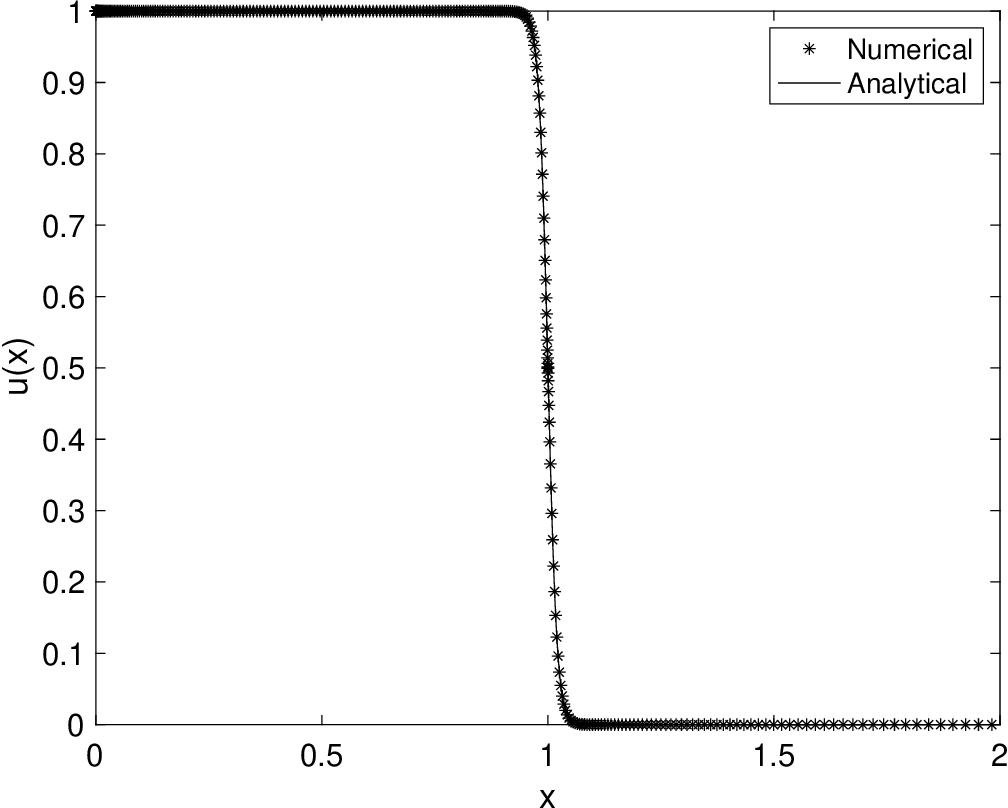}
\includegraphics[height=2.0in]{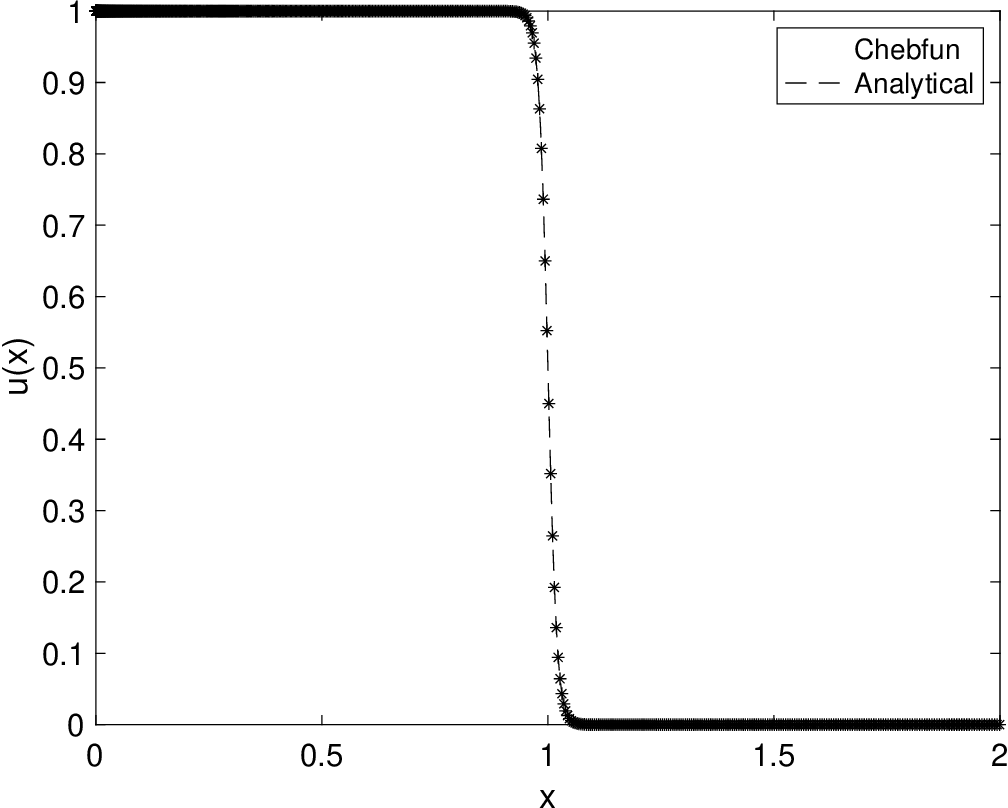}
\caption{Comparison of the error in approximating the function for a solution using \textit{chebfun} with the results of the one-point overlap method.}
\label{1point_comparison}
\end{figure}

Figure \ref{1point_comparison}(a) shows the performance of the 1-point algorithm for solving the above BVP for $\theta = 0.02$ and N=200 points. Using a single global mapping to cluster points n the vicinity of y=1 required 6500 points to achieve an acceptable level of normalized error. In Fig. \ref{1point_comparison}(b) we show the results when using \textit{chebfun} with splitting enabled. The code converges quickly, but it can be observed that the distribution of points is somewhat different, with fewer points in the region of width $\theta$ around y=1. This is a result of the automatic grid bisection and regrouping feature in \textit{chebfun}. 

\begin{figure}[h]
\centering
\includegraphics[width=0.49\textwidth]{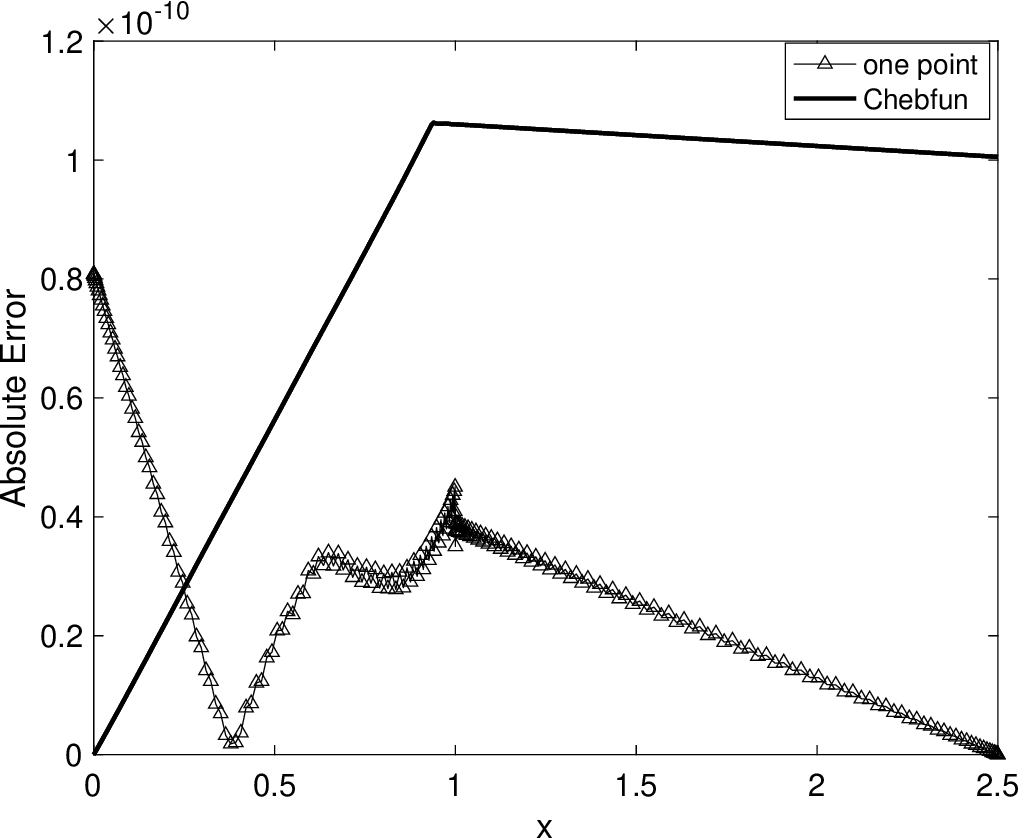}
\includegraphics[width=0.49\textwidth]{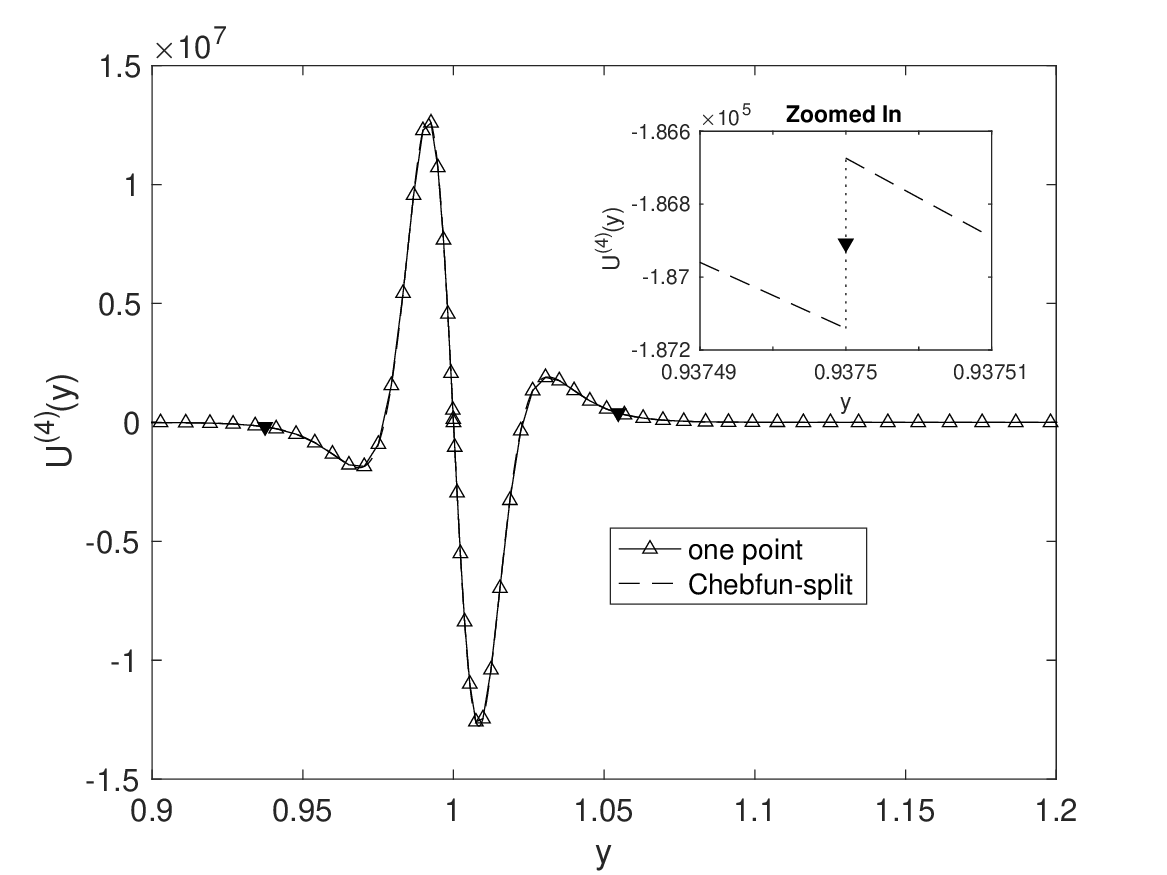}
\caption{Comparison of the error in approximating the function for solution using \textit{chebfun} with the results of the one-point overlap method.}
\label{1point_comparison_error}
\end{figure}
Figure \ref{1point_comparison_error}(a) shows the error between the solutions computed from the one-point overlap method and \textit{chebfun} with splitting. The errors are comparable, with the maximum error in the 1-point method being highest near y=1 and decreasing thereafter. \text{chebfun} on the other hand has peak error away from the gradient region due to the choice of locations of automatic interval bisection. 

Figure \ref{1point_comparison_error}(b) compares the value of the fourth derivatives of the function when calculated using the two methods. It can be seen that while the present method preserves the continuity of this higher order derivative due to its inherent formulation, the solution from \textit{chebfun} has a jump at subdomain boundaries (see inset). This is one of the main advantages of using the algorithm when compaed to chebfun with splitting. The present algorithm uses more points (300  compared to 185) but preserves the full extent of higher order representation afforded by Chebyshev collocation. At the same time, it uses vastly fewer points (300 vs 6500) for solution that does not use domain decomposition and relies on chebdun without splitting enabled.  

\subsection{Algorithm 2: Two point overlap}
Intuitively, one might consider whether we can increase the number of overlapping points to improve the convergence and accuracy beyond that offered by algorithm 1. In the 1-point algorithm and the piecewise continuous approximation employed by \textit{chebfun} with splitting the boundary overlap point is the sole point through which information about the two separate polynomial representations of the function in the subdomains is exchanged, leading to a single row with completely filled elements. Indeed, the unique properties of Chebyshev collocation methods arise from the filled operator matrices which ensure fast convergence with no iteration required unlike finite difference schemes, which use information from a stencil of nearby points. More points from each subdomain that overlap will incorporate information from each subdomain and would be expected to improve convergence.

Once a choice of mapping is made for each subdomain, we can construct the derivative matrix as described in algorithm 1. The simplest way to ensure the existence of an interval containing two points from either subdomain is to choose mappings that ensure two points from one subdomain coincide with two points from the other subdomain. We can satisfy the above requirements if the nonlinear mappings are the same in the two domains; however other choices are possible and we illustrate one such case. Suppose we map the subdomains using the mappings $y_1$ and $y_2$  as follows
\begin{equation}
\begin{array}{l}
y_{1}=\frac{H_{1}}{2}(1-x) \\
y_{2}=\frac{H_{2}}{2}(1-x) + 
\end{array}
\end{equation}
We have the derivative matrix $D_1 = -\frac{2}{H_1}D$.  Linear mappings allow for stretching and scaling in a convenient manner so that points from either subdomain can be made to coincide, i.e. the last two points of subdomain $H_1$ coincide with the first two points of subdomain $H_2$. The mapping on $y_2$ can be obtained by shifting the mapping on $y_1$ to the right with a step $H_1$. The global matrix can be written as:  
\begin{equation}
D=\left[\begin{array}{ll}
\frac{-2D}{H_{1}} & \\
& \frac{-2D}{H_{1}}
\end{array}\right]
\end{equation}
To show how to assemble the global matrix, consider an example with four points in each domain as follows.  
\begin{equation}
\left[\begin{array}{llllllll}
a_{11} & a_{12} & a_{13} & a_{14} & & \\
a_{21} & a_{22} & a_{23} & a_{24} & & \\
a_{31} & a_{32} & a_{33} & a_{34} & & \\
a_{41} & a_{42} & a_{43} & a_{44} & & \\
& & & &  b_{11} & b_{12}  & b_{13}  & b_{14} \\
& & & & b_{21} & b_{22}  & b_{23}  & b_{24} \\
& & & & b_{31} & b_{32}  & b_{33}  & b_{34} \\
& & & & b_{41} & b_{42}  & b_{43}  & b_{44}
\end{array}\right]\left[\begin{array}{c}
u_{1} \\
u_{2} \\
u_{3} \\
u_{4} \\
v_{1} \\
v_{2} \\
v_{3} \\
v_{4}
\end{array}\right]=\left[\begin{array}{l}
d u_{1} \\
d u_{2} \\
d u_{3} \\
d u_{4} \\
d v_{1} \\
d v_{2} \\
d v_{3} \\
d v_{4}
\end{array}\right]
\end{equation}
The locations corresponding to the points $u_3$, $u_4$, $v_1$, and $v_2$ are the overlapping grid points. In other words, $u_3=v_1$ and $u_4= v_2$. Two rows need to be removed from the matrix to ensure that it is non-singular. As in algorithm 1, the extra equations arise from separate enforcement of the governing equations using two different polynomial representations corresponding to the two subdomains. These extra equations are removed by multiplying them by a constant weighting factor and adding them to weighted equations from the other subdomain. In following this algorithm, we shift the matrix D to up and left simultaneously by two elements and use half of the value from both regions for the overlapping points. As before, we note that weights other than 0.5 for each side are possible. This results in the following global matrix, where the averaged values for the overlapping points can be noted:
\begin{equation}
\left[\begin{array}{llllll}
a_{11} & a_{12} & a_{13} & a_{14} & & \\
a_{21} & a_{22} & a_{23} & a_{24} & & \\
a_{31}/2 & a_{32}/2 & \frac{b_{13} +  b_{11}}{2}  &\frac{b_{14} + b_{12}}{2}& & \\
a_{41}/2 & a_{42}/2 & \frac{b_{23} +  b_{21}}{2}  &\frac{b_{24} +  b_{22}}{2}
& & \\
& & b_{31} & b_{32}  & b_{33}  & b_{34} \\
& & b_{41} & b_{42}  & b_{43}  & b_{44}
\end{array}\right]\left[\begin{array}{c}
u_{1} \\
u_{2} \\
u_{3} \\
u_{4} \\
v_{3} \\
v_{4}
\end{array}\right]=\left[\begin{array}{l}
d u_{1} \\
d u_{2} \\
d u_{3} \\
d u_{4} \\
d v_{3} \\
d v_{4}
\end{array}\right]
\label{2pt_m}
\end{equation}
The second-order derivative matrix can be obtained via $D_2^* = D^*D^*$,  where $D^*$ is the derivative operator defined above. While the example above used two linear mappings, one can *(in principle) choose to use nonlinear mappings for both regions or combine a linear mapping with another nonlinear mapping. In practice, however, trying to make points from each subdomain coincide imposes severe restrictions on the type of mappings and may not be practical in all situations. We develop an approach to relax this requirement in the next section.

\subsubsection{Example: Eigenvalue problem of hydrodynamic instability for the weakly diffusive interface between two fluid streams of different viscosity }
An example of the performance of the two-point algorithm is provided below. We consider the laminar flow of two fluids of different viscosity and equal density in a round tube. Fluid 1 occupies the core region of diameter $2r_i$ while fluid 2 flows in the annular space $r_i < r < 1$. For immiscible fluids, an interface of zero thickness can be defined across which conditions for continuity of velocity and stresses can be defined. This situation strongly suggests domain decomposition as the solution approach, enabling clustering points near the interface in the two fluid subdomains through separate mappings. However, when the fluids are miscible, weak diffusion leads to a thin diffusive interface of extent $\delta$ where the viscosity changes sharply from the inner to the outer fluid value. This leads to strong but continuous changes in the velocity profile in the region of the diffusive interface, which are known to be more challenging to resolve than the zero-thickness situation. 
Hydrodynamic instabilities that arise in this configuration are modeled as traveling wave packets with spatially varying amplitude. The problem can be formulated as an eigenvalue problem. The perturbed Navier-Stokes equations can be written as 

$$
\begin{aligned}
& \frac{\mathrm{d} \hat{v}_r}{\mathrm{~d} r}+\frac{\hat{v}_r}{r}+\frac{\beta \hat{v}_\theta}{r}+\alpha \hat{v}_z=0, \\
& \operatorname{Re}\left[-\omega \hat{v}_r+\alpha \bar{v}_z \hat{v}_r\right]=\frac{\mathrm{d} \hat{p}}{\mathrm{~d} r}-\mathrm{ie}^{M \bar{c}}\left[\frac{\mathrm{d}^2 \hat{v}_r}{\mathrm{~d} r^2}+\frac{1}{r} \frac{\mathrm{d} \hat{v}_r}{\mathrm{~d} r}-\left(\frac{\beta^2+1}{r^2}+\alpha^2\right) \hat{v}_r\right. \\
& \left.-\frac{2 \beta}{r^2} \hat{v}_\theta+2 M \frac{\mathrm{d} \bar{c}}{\mathrm{~d} r} \frac{\mathrm{d} \hat{v}_r}{\mathrm{~d} r}+M \alpha \frac{\mathrm{d} \bar{v}_z}{\mathrm{~d} r} \hat{c}\right], \\
& \operatorname{Re}\left[-\omega \hat{v}_\theta+\alpha \bar{v}_z \hat{v}_\theta\right]=\frac{-\beta \hat{p}}{r}-\mathrm{ie}^{M \bar{c}}\left[\frac{\mathrm{d}^2 \hat{v}_\theta}{\mathrm{d} r^2}+\frac{1}{r} \frac{\mathrm{d} \hat{v}_\theta}{\mathrm{d} r}-\left(\frac{\beta^2+1}{r^2}+\alpha^2\right) \hat{v}_\theta\right. \\
& \left.-\frac{2 \beta}{r^2} \hat{v}_r+M \frac{\mathrm{d} \bar{c}}{\mathrm{~d} r}\left(\frac{\mathrm{d} \hat{v}_\theta}{\mathrm{d} r}-\frac{\hat{v}_\theta}{r}-\frac{\beta \hat{v}_r}{r}\right)\right], \\
& \operatorname{Re}\left[-\omega \hat{v}_z+\alpha \bar{v}_z \hat{v}_z+\frac{\mathrm{d} \bar{v}_z}{\mathrm{~d} r} \hat{v}_r\right]=-\alpha \hat{p}-\mathrm{ie}^{M \bar{c}}\left[\frac{\mathrm{d}^2 \hat{v}_z}{\mathrm{~d} r^2}+\frac{1}{r} \frac{\mathrm{d} \hat{v}_z}{\mathrm{~d} r}-\left(\frac{\beta^2}{r^2}+\alpha^2\right) \hat{v}_z\right. \\
& +M \frac{\mathrm{d} \bar{c}}{\mathrm{~d} r}\left(\frac{\mathrm{d} \hat{v}_z}{\mathrm{~d} r}-\alpha \hat{v}_r\right)+M \frac{\mathrm{d} \bar{v}_z}{\mathrm{~d} r} \frac{\mathrm{d} \hat{c}}{\mathrm{~d} r} \\
& \left.+M \hat{c}\left(\frac{\mathrm{d}^2 \bar{v}_z}{\mathrm{~d} r^2}+\frac{1}{r} \frac{\mathrm{d} \bar{v}_z}{\mathrm{~d} r}+M \frac{\mathrm{d} \bar{c}}{\mathrm{~d} r} \frac{\mathrm{d} \bar{v}_z}{\mathrm{~d} r}\right)\right], \\
& P e\left[-\omega \hat{c}+\alpha \bar{v}_z \hat{c}+\frac{\mathrm{d} \bar{c}}{\mathrm{~d} r} \hat{v}_r\right]=-\mathrm{i}\left[\frac{\mathrm{d}^2 \hat{c}}{\mathrm{~d} r^2}+\frac{1}{r} \frac{\mathrm{d} \hat{c}}{\mathrm{~d} r}-\left(\frac{\beta^2}{r^2}+\alpha^2\right) \hat{c}\right] \text {. } \\
&
\end{aligned}
$$

Here, $\hat{v}_r(r)$, $\hat{v}_z(r)$ and $\hat{v}_\theta (r)$ are the velocity perturbations in the radial, axial and azimuthal directions respectively, while $\hat{p}(r)$ and $\hat{c}(r)$ are the perturbations in the pressure and species concentration. The base state profiles for unperturbed quantities are given by $\bar{v}_z(r)$, $\bar{p}(r)$ and $\bar{p}(r)$, with $\bar{v}_r=\hat{v}_\theta =0$. The quantities $\alpha$ and $\omega$ are the complex wavenumber and frequency of traveling waves, while $Re$, $Pe$ and $M$ are specified input parameters (Reynolds number, Peclet number, and viscosity ratio) that, along with the unperturbed profiles, govern the relationship between the wavenumber $\alpha$ and frequency $\omega$.    

The boundary conditions at the centerline and the pipe radius ($r=1$) are are as follows. For both cases, at the center-line, we consider the single-valueness of velocity, together with continuity to derive the center-line conditions for different values of azimuthal wavenumber $\beta$  (\cite{khorrami1989application}):

\begin{equation}
\begin{aligned}
& \beta=0: \quad \frac{\mathrm{d} \hat{v}_z}{\mathrm{~d} r}=0, \quad \hat{v}_r=0, \quad \hat{v}_\theta=0, \quad \frac{\mathrm{d} \hat{p}}{\mathrm{~d} r}=0, \quad \frac{\mathrm{d} \hat{c}}{\mathrm{~d} r}=0, \\
& \beta=1: \quad \hat{v}_z=0, \quad \hat{v}_r+\hat{v}_\theta=0, \quad 2 \frac{\mathrm{d} \hat{v}_r}{\mathrm{~d} r}+\frac{\mathrm{d} \hat{v}_\theta}{\mathrm{d} r}=0, \quad \hat{p}=0, \quad \hat{c}=0, \\
& \beta 2: \quad \hat{v}_z=0, \quad \hat{v}_r=0, \quad \hat{v}_\theta=0, \quad \hat{p}=0, \quad \hat{c}=0 .
\end{aligned}
\end{equation}

Together, the above equations constitute an eigenvalue problem, with the angular frequency representing the eigenvalue while the velocity and concentration disturbances are the eigenfunctions. In other words, these equations represent the dispersion relation: 
\begin{equation}
    D(\omega, k, Re, M, Sc, \theta , \theta _\mu , \delta ) =0
\end{equation}

The eigenvalues  $\omega$ are the complex-valued frequencies of traveling wave packets with wavenumber $k$. The eigenvectors for each of the eigenmodes can be summed to yield the perturbation velocities ($\hat{v}_z$, $\hat{v}_r$ and $\hat{v}_\theta$). One proposed set of base profiles for the unperturbed axial velocity and species concentration is:
\begin{equation}
\bar{c}(r)=0.5+0.5 \operatorname{erf}\left(\frac{r-a}{\delta}\right).
\end{equation}
The axial base flow must be evaluated numerically from the dimensionless axial momentum equation in the $z$-direction
\begin{equation}
\frac{\mathrm{d} \bar{p}}{\mathrm{~d} z}=\mathrm{e}^{M \bar{c}}\left[\frac{\mathrm{d}^2 \bar{v}_z}{\mathrm{~d} r^2}+\frac{1}{r} \frac{\mathrm{d} \bar{v}_z}{\mathrm{~d} r}+M \frac{\mathrm{d} \bar{v}_z}{\mathrm{~d} r} \frac{\mathrm{d} \bar{c}}{\mathrm{~d} r}\right].
\end{equation}
where assume the viscosity $\mu$ to be an exponential function of the concentration 
\begin{equation}
\mu=\mu_2 \mathrm{e}^{M c}, \quad M=\ln \frac{\mu_1}{\mu_2} .
\end{equation}

\begin{figure}
\centering
\includegraphics[height=1.8in]{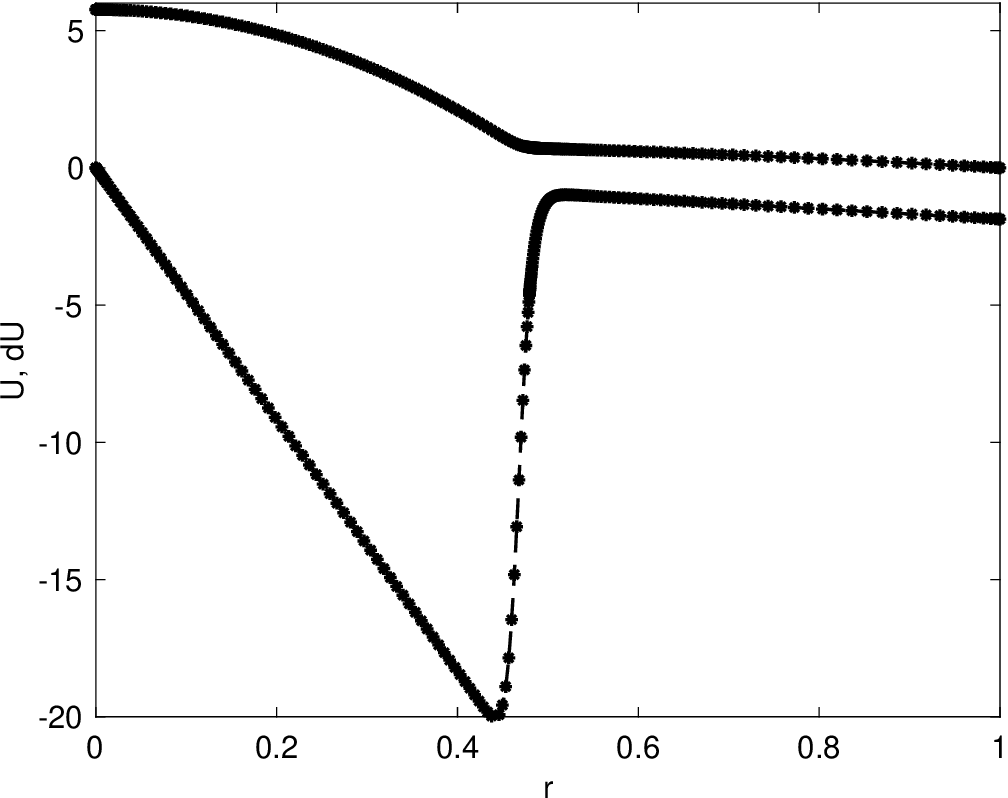}
\includegraphics[height=1.8in]{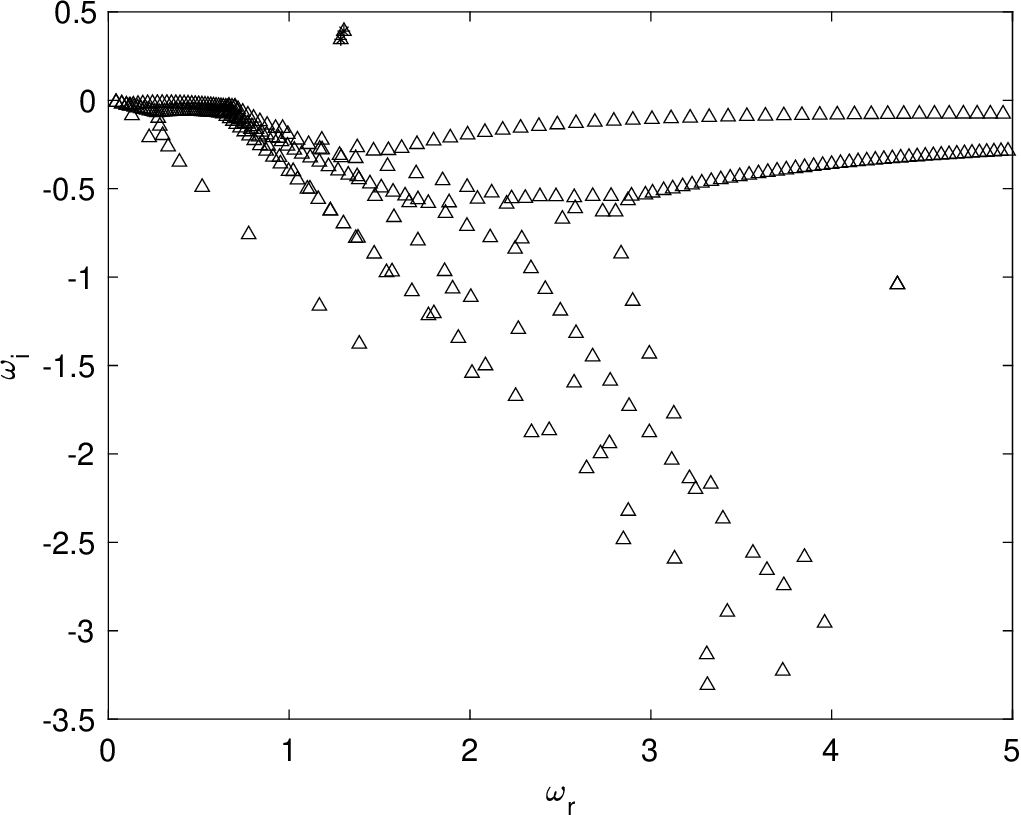}
\caption{3 D channel flow instability with $\delta = 0.02$: (a) velocity profile and its derivative (b) Spectral distribution when wave number $k=1-i$}
\label{fig:Selvam_comp}
\end{figure}

Figure \ref{fig:Selvam_comp} shows the profiles of velocity, viscosity and the velocity derivative for $Re=$, $Pe=$, and $\delta = 0.02$. This problem has been previously solved by Selvam et al. \cite{Selvam2009} using a global mapping (i.e. no domain decomposition) for a range of flow conditions. The 'neutral' curve, i.e. the values of Re for which the temporal growth rate, i.e. Im($\omega$) =0 is shown in Fig. \ref{fig:Selvam_comp}(b), along with calculations based on the two-point overlap algorithm. It is clear that the algorithm performs at least as well as current global methods to capture such weak discontinuities in the interior of the domain. We now show that this extends to even sharper discontinuities corresponding to smaller values of $\delta$. 

\begin{figure}
\centering
\includegraphics[height=2.3in]{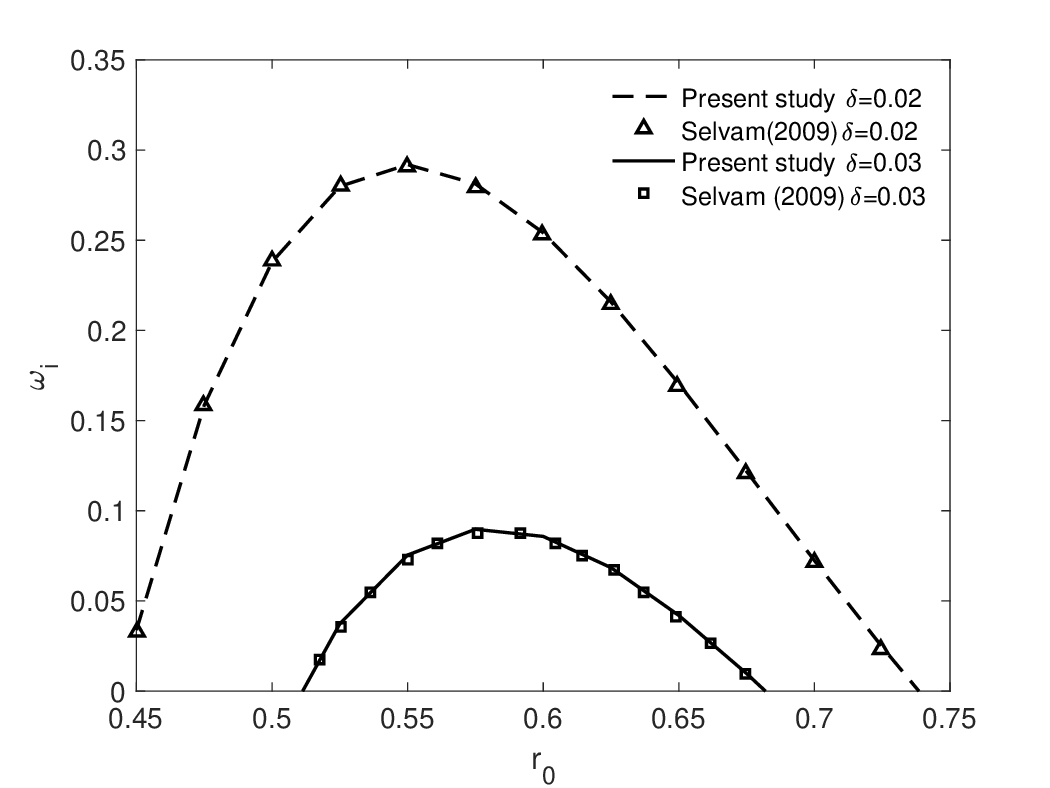}
\includegraphics[height=2.3in]{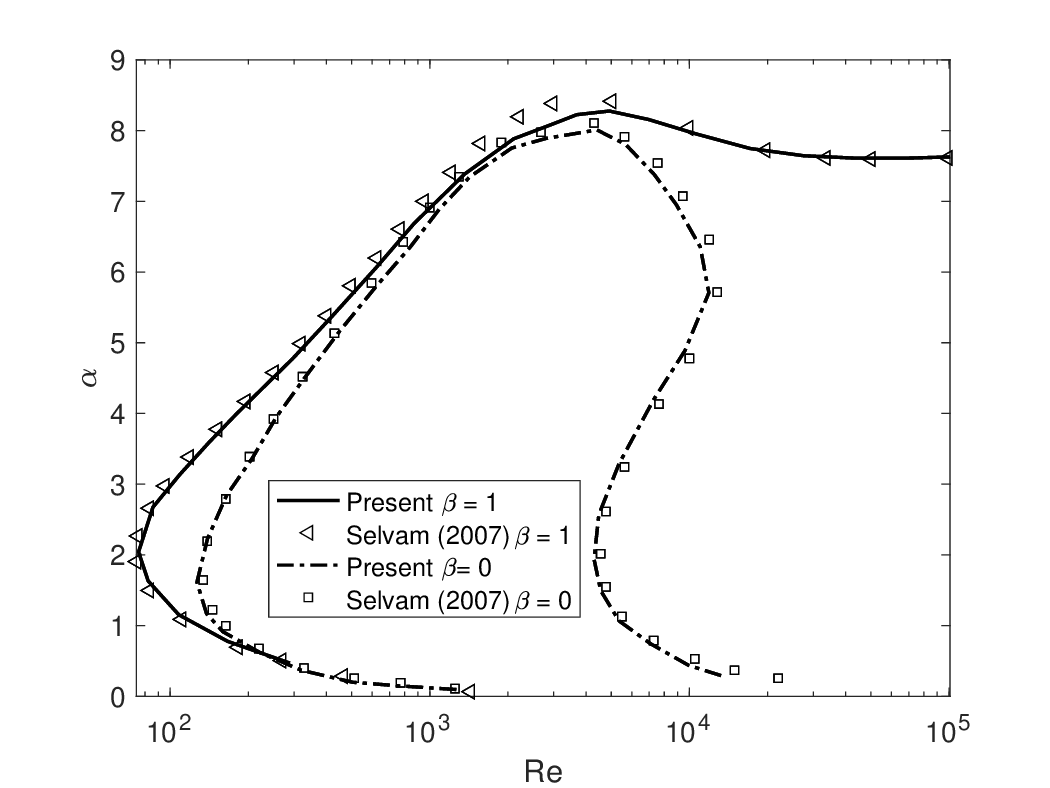}
\caption{(a) Comparison of our calculations with the core-annular pipe flow calculations in Fig.3(a) of \cite{Selvam2009} for the axisymmetric mode. The abscissa $r_0$ corresponds to various locations of the diffusive interface of width $\delta$, for the parameter matrix $(M,Re,Sc,\delta^*)=(25,48,7500,[0.02,0.03])$. (b) Comparison with the results from Selvam in Fig.7(2007), considering both helical and axisymmetric modes when the ratio of annular to core fluid viscosity is 2.718.}
\label{M25_Validation}
\end{figure}

\clearpage
\subsection{Algorithm 3: Pseudo Multi points overlapping method}
The previous two sections have shown that while the accuracy of higher-order collocation methods can be maintained by having overlapping points, the approach comes with some constraints. Namely, the overlapping points correspond to subdomain boundaries, and matching two points from one domain on points of another is convenient only with linear mappings. Now we consider how these constraints may be relaxed, i.e. we consider how the overlapping points might lie in the interiors of the two subdomains, while also allowing for nonlinear mappings. 
We continue to consider two regions $H_1$ and $H_2$, while the two mappings from the calculation domain into $[-1,1]$ are 
\begin{equation}
\begin{array}{l}
y_{1}=\frac{H_{1}(1-x)}{2} \\
y_{2}=\frac{L(1-x)}{b+x} + H_1 - \delta
\end{array}
\end{equation}
where $L$ and $b$ are two free parameters used to control the grid point distribution in the region $[0, H_1]$. $\delta$ is the region of overlap between the two regions, with no restriction on the distribution of points in this region. For illustration, we consider a situation with five points in each subdomain, with a total of three points combined in the overlap region. The process of building the derivative matrix is illustrated below. 

When considering the two subdomains separately, The symbolical representation of the numerical derivative is written as follows
\begin{equation}
\left[\begin{array}{llllllllll}
a_{11} & a_{12} & a_{13} & a_{14} & a_{15} & & \\
a_{21} & a_{22} & a_{23} & a_{24} & a_{25}& & \\
a_{31} & a_{32} & a_{33} & a_{34} & a_{35}& & \\
a_{41} & a_{42} & a_{43} & a_{44} & a_{45}& & \\
a_{41} & a_{42} & a_{43} & a_{44} & a_{55}& & \\
& & & & & b_{11} & b_{12}  & b_{13}  & b_{14} & a_{15}\\
& & & & & b_{21} & b_{22}  & b_{23}  & b_{24} & a_{25}\\
& & & & & b_{31} & b_{32}  & b_{33}  & b_{34} & a_{35}\\
& & & & & b_{41} & b_{42}  & b_{43}  & b_{44} & a_{45}\\
& & & & & b_{51} & b_{52}  & b_{53}  & b_{54} & a_{55}
\end{array}\right]\left[\begin{array}{c}
u_{1} \\
u_{2} \\
u_{3} \\
u_{4} \\
u_{5} \\
v_{1} \\
v_{2} \\
v_{3} \\
v_{4} \\
v_{5}
\end{array}\right]=\left[\begin{array}{l}
d u_{1} \\
d u_{2} \\
d u_{3} \\
d u_{4} \\
d u_{5} \\
d v_{1} \\
d v_{2} \\
d v_{3} \\
d v_{4} \\
d v_{5}
\end{array}\right]
\label{3pt_D}
\end{equation}
To account for the three points in the overlap region, \Cref{3pt_D} can be rewritten as 
\begin{equation}
\left[\begin{array}{llllllllll}
a_{11} & a_{12} & 0 & a_{13} & 0 & a_{14} & 0 & a_{15} & &\\
a_{21} & a_{22} & 0 & a_{23} & 0 & a_{24} & 0 & a_{25} & &\\
& & b_{11} & 0 & b_{12} & 0  & b_{13} & 0  & b_{14} & b_{15}\\
a_{31} & a_{32} & 0 & a_{33} & 0 & a_{34} & 0 & a_{35} & &\\
& & b_{21} & 0 & b_{22} & 0  & b_{23} & 0  & b_{24} & b_{25}\\
a_{41} & a_{42} & 0 & a_{43} & 0 & a_{44} & 0 & a_{45} & &\\
& & b_{31} & 0 & b_{32} & 0  & b_{33} & 0  & b_{34} & b_{35}\\
a_{51} & a_{52} & 0 & a_{53} & 0 & a_{54} & 0 & a_{55} & &\\
& & b_{41} & 0 & b_{42} & 0  & b_{43} & 0  & b_{44} & b_{45}\\
& & b_{51} & 0 & b_{52} & 0  & b_{53} & 0  & b_{54} & b_{55}
\end{array}\right]\left[\begin{array}{c}
u_{1} \\
u_{2} \\
v_{1} \\
u_{3} \\
v_{2} \\
u_{4} \\
v_{3} \\
u_{5} \\
v_{4} \\
v_{5}
\end{array}\right]=\left[\begin{array}{l}
du_{1} \\
du_{2} \\
dv_{1} \\
du_{3} \\
dv_{2} \\
du_{4} \\
dv_{3} \\
du_{5} \\
dv_{4} \\
dv_{5}
\end{array}\right]
\label{3pt_D1}
\end{equation}

\begin{figure}
\centering
\includegraphics[height=2.5in]{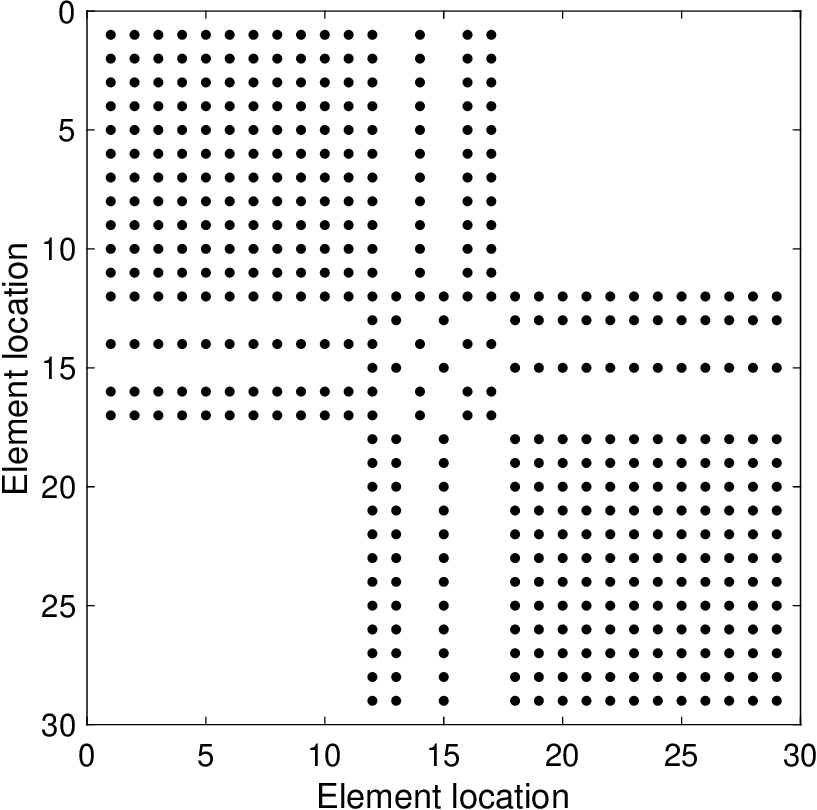}
\includegraphics[height=2.5in]{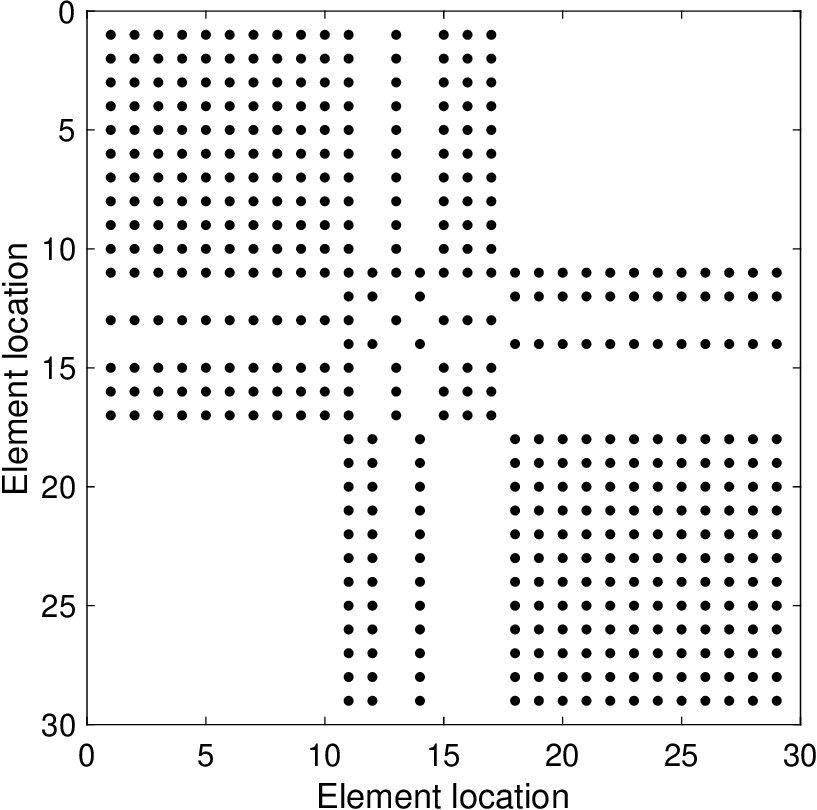}
\includegraphics[height=2.5in]{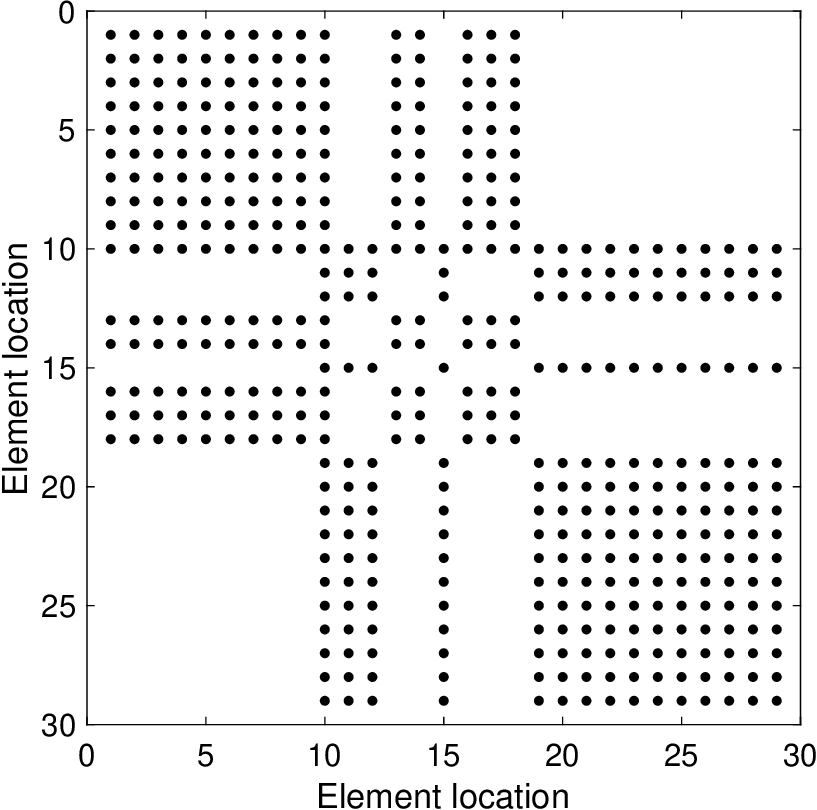}
\includegraphics[height=2.5in]{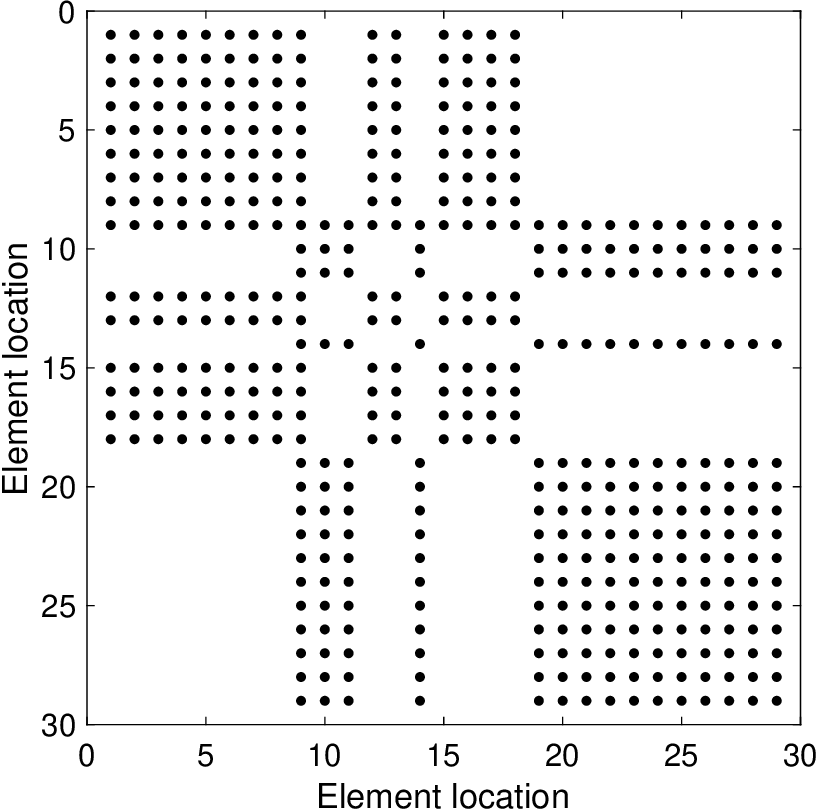}
\caption{Non-zero matrix elements distribution of the multi-points overlapping method with a total 15 points in each domain: (a) three points overlapping global matrix; (b) four points overlapping global matrix;  (c) five points overlapping global matrix; (d) six points overlapping global matrix.}
\label{mt_pt_fig}
\end{figure}

To avoid the case that the global matrix in \Cref{3pt_D1} is a singular matrix, we can assume the grid points corresponding to $u_2$ and $v_1$ are coincident. We can also make a similar assumption for grid points $u_5$ and $v_4$. We show the distribution of non-zero matrix elements  for this `pseudo-multi-point overlap' method in \Cref{mt_pt_fig}.  This term is used to refer to the fact that, though there are multiple points in the overlap region, only one grid point belongs to both domains, as evident from point eqn. \cref{3pt_D1}. The corresponding code to provide details  of the derivative matrix assembly process is given in the Appendix.

\Cref{mt_pt_fig} shows the distribution of non-zero matrix elements for the pseudo-multipoint overlapping method with 15 points in each domain. In \Cref{mt_pt_fig}(a), there are three points in the overlapping region, which can be seen from between the 12th and 20th nodes. In \Cref{mt_pt_fig}(b,c,d) there are four, five, and six point in the overlapping region, respectively. 

\begin{figure}
\centering
\includegraphics[height=2in]{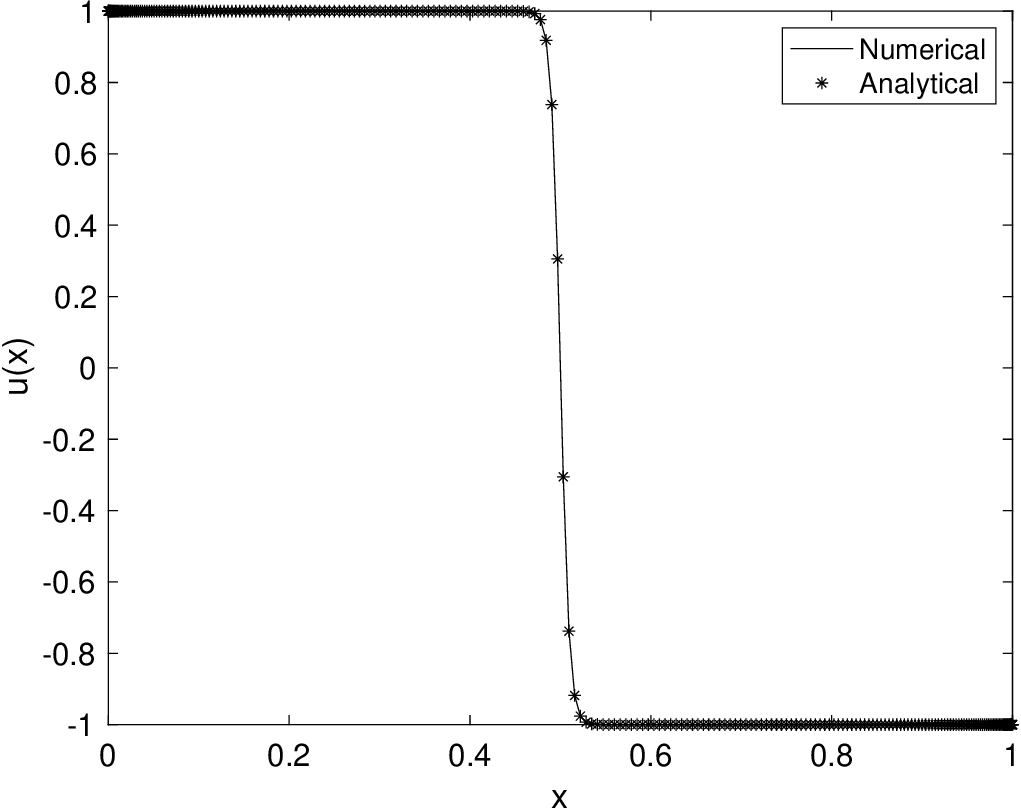}
\includegraphics[height=2in]{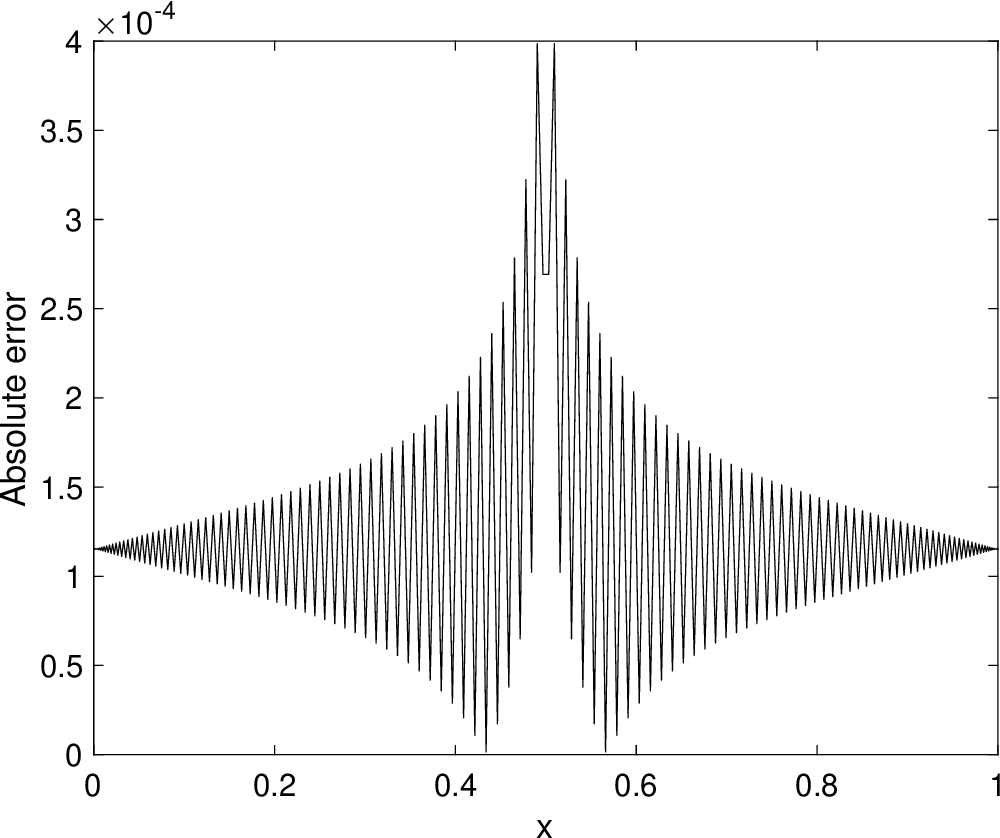}
\caption{Absolute error and solution for the Burgers equation via global Chebyshev method with nodes $N = 300$ and $\nu = 5 \times 10^{-3}$.}
\label{fig:Burgers_global}
\end{figure}

Figure \ref{fig:Burgers_global}(a)  shows the solution to the analytical solution to the Burgers equation as well as the solution computed from a fully global method implementing a single mapping given by 

\begin{figure}
\centering
\includegraphics[height=2in]{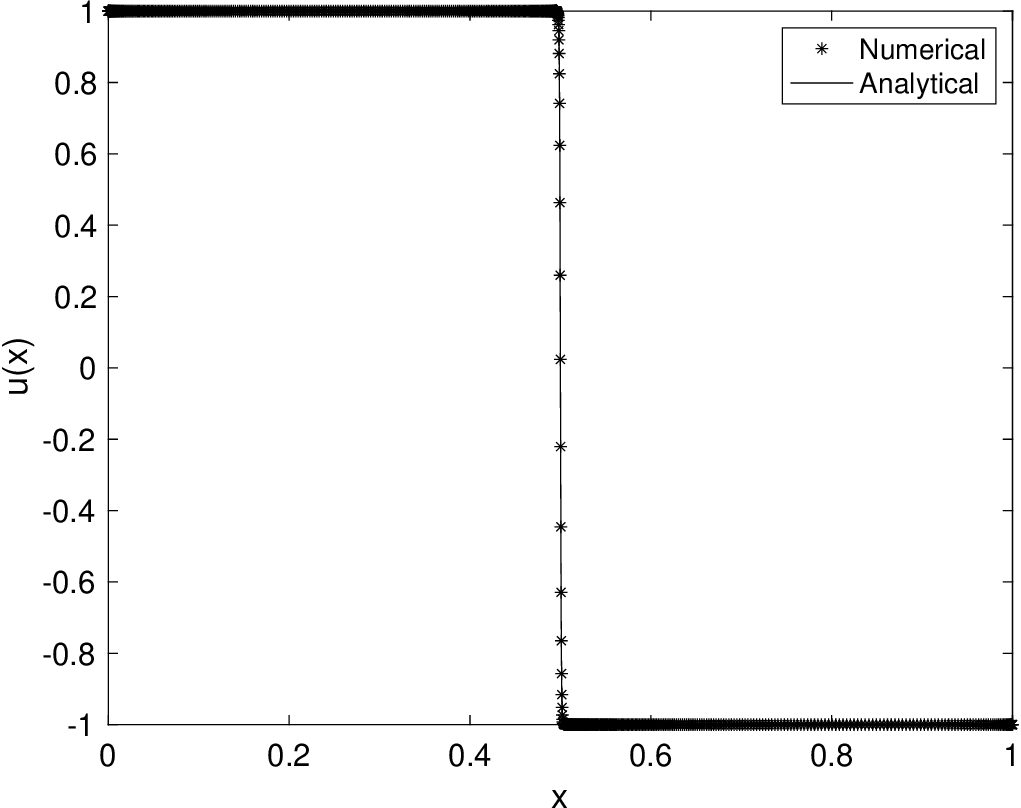}
\includegraphics[height=2in]{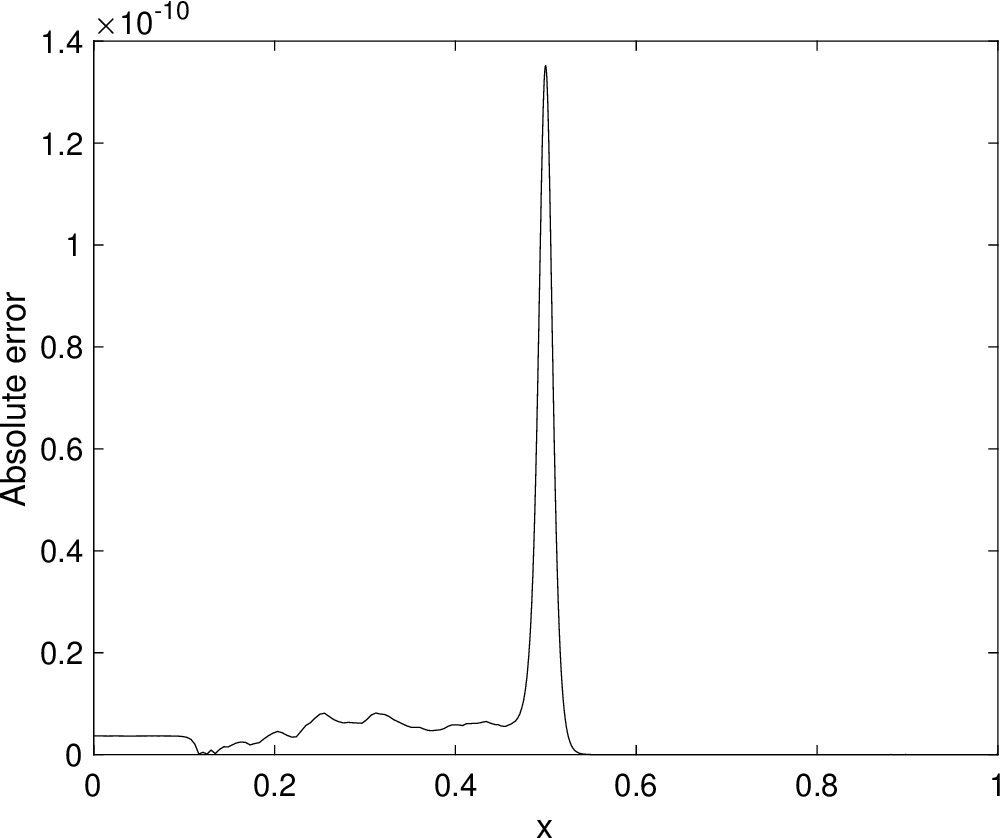}
\caption{Absolute error and solution for the Burgers equation via semi-global Chebyshev method with nodes $N = 299$  and $\nu = 5 \times 10^{-3}$ for two points overlapping.}
\label{fig:Burgers_2pt}
\end{figure}

While N=300 points are used over the calculation domain, it is clear that only a few points can be located in the region of the sharpest gradient. The error between the numerical and analytical solutions is shown in Fig. \ref{fig:Burgers_global}(b) and suggests a peak error of $10^{-4}$.  Figure \ref{fig:Burgers_2pt}(a) and (b) show the calculated solution and the error from the analytical solution, when using the pseudo-multipoint overlap method, for the case of two points in the overlap interval, with 150 points in each subdomain. \textbf{The mappings used for the subdomains are the same as in  *****}. It is evident that there are more points in the gradient region, and further, the error has gone down by several orders of magnitude. Fig. \ref{fig:Burgers_global} shows that when additional points are included in the overlap region (five in this example), the peak error goes down further by one order of magnitude. 

\begin{figure}
\centering
\includegraphics[height=2in]{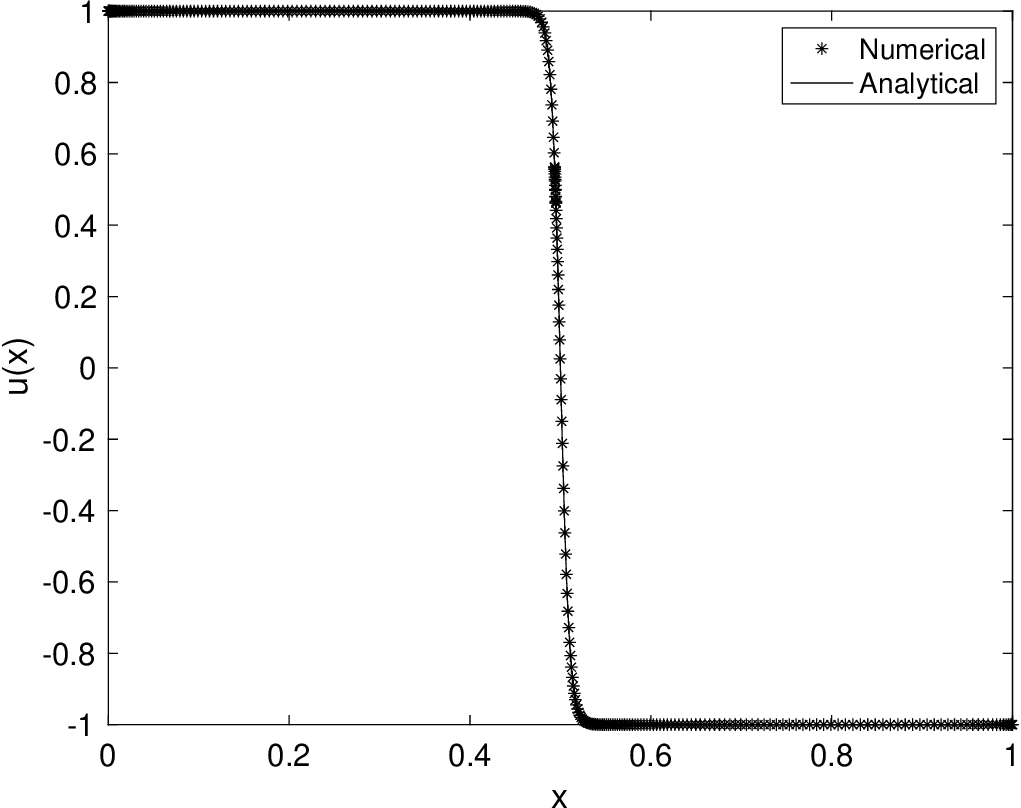}
\includegraphics[height=2in]{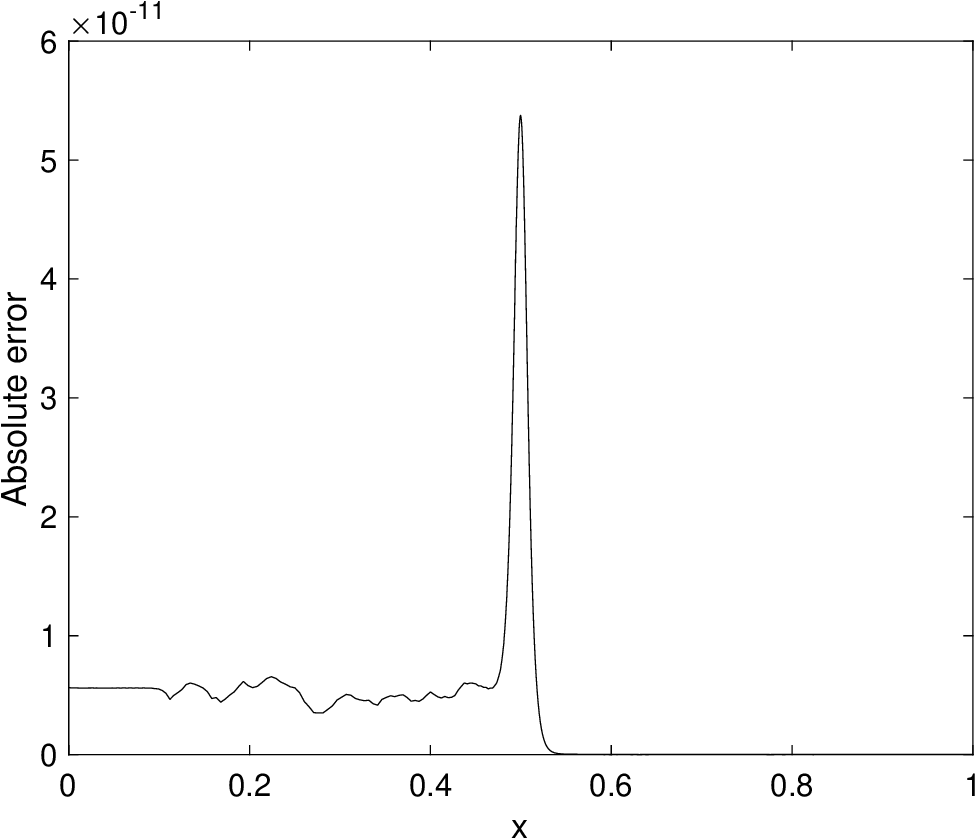}
\caption{Absolute error and solution for the Burgers equation via semi-global Chebyshev method with nodes $N = 299$  and $\nu = 5 \times 10^{-3}$ for five points overlapping.}
\label{fig:Burgers_5pt}
\end{figure}

\begin{figure}
\centering
\includegraphics[height=2in]{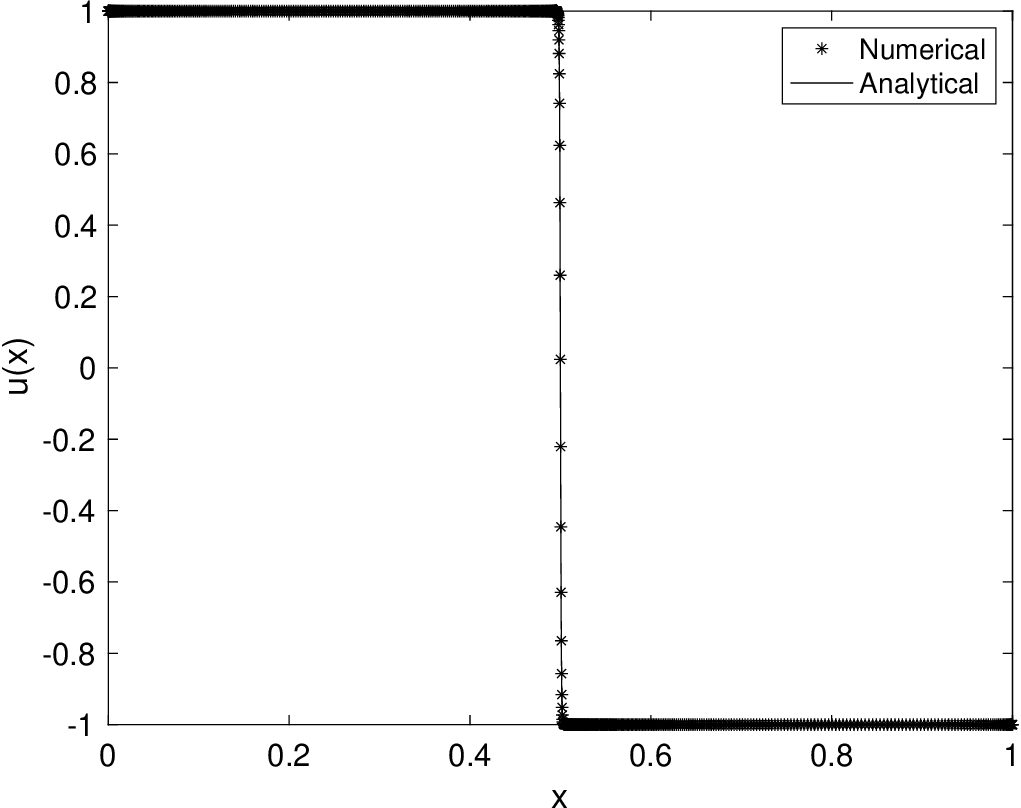}
\includegraphics[height=2in]{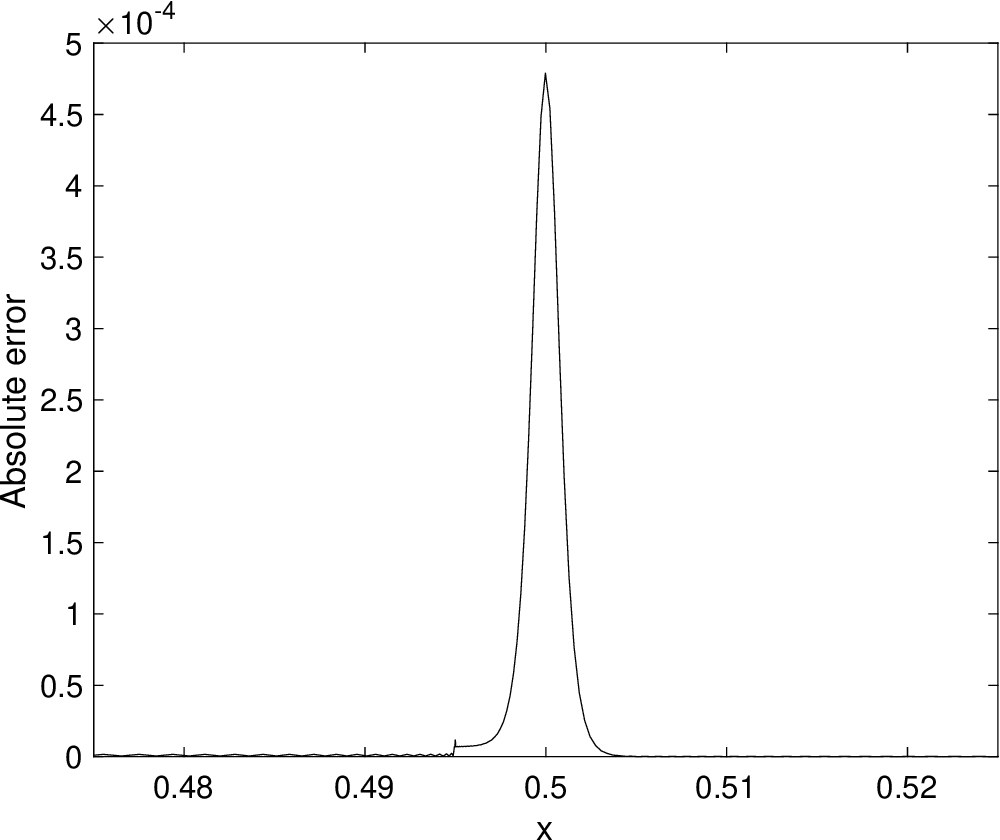}
\caption{Absolute error and solution for the Burgers equation via five points overlapped Chebyshev method with nodes $N = 299$  and \textbf{$\nu = 5 \times 10^{-4}$} for five points overlapping.}
\end{figure}

\subsubsection{Case Study: Nonlinear Boundary Value Problem:} 
We solve the stationary Burgers equation on the interval $[0,1]$ with Robin boundary conditions \cite{reyna1995exponentially} :
$$
\begin{aligned}
& \nu u^{\prime \prime}(x)-u(x) u^{\prime}(x)=0, \\
& \nu u^{\prime}(0)-\kappa(u(0)-\alpha)=0, \\
& \nu u^{\prime}(1)+\kappa(u(1)+\alpha)=0,
\end{aligned}
$$

which has a nontrivial solution
$$
u(x)=-\beta \tanh \left(\frac{1}{2} \beta \nu^{-1}\left(x-\frac{1}{2}\right)\right),
$$
where $\beta$ satisfies
$$
-\frac{1}{2} \beta^2 \operatorname{sech}^2\left(\frac{1}{4} \beta \nu^{-1}\right)+\kappa\left[\alpha-\beta \tanh \left(\frac{1}{4} \beta \nu^{-1}\right)\right]=0 .
$$

We choose $\nu=5 \times 10^{-3}, \alpha=1$, and $\kappa=2$.

\begin{figure}
\centering
\includegraphics[height=2.0in]{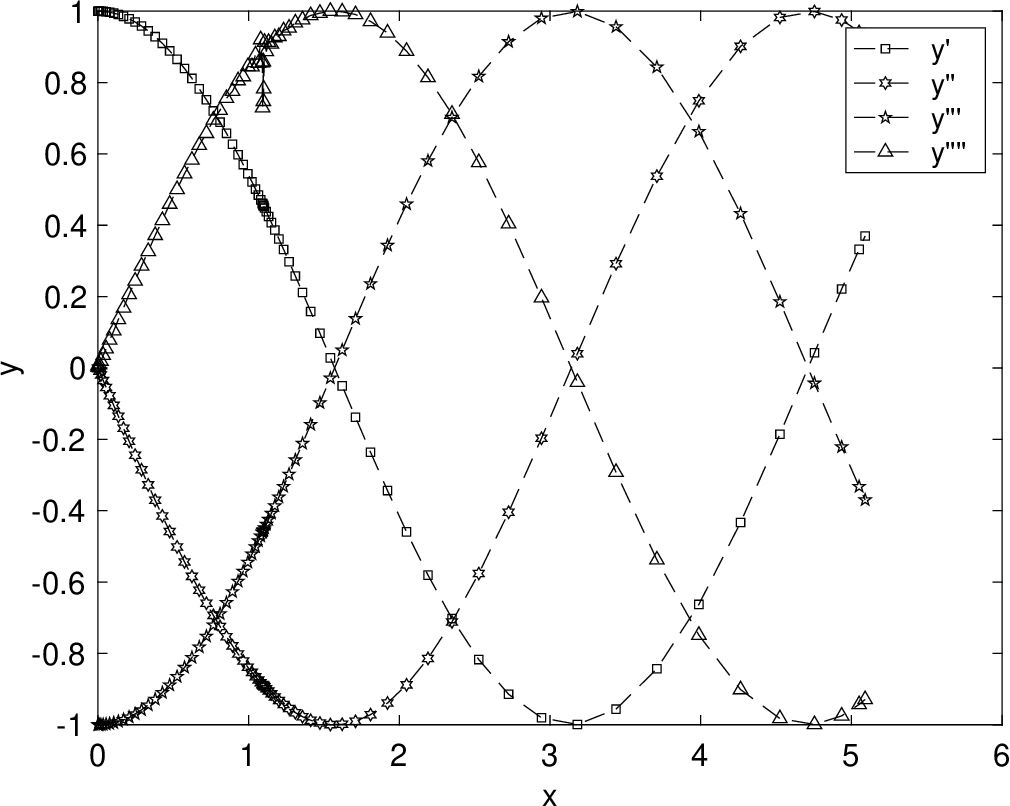}
\includegraphics[height=2.0in]{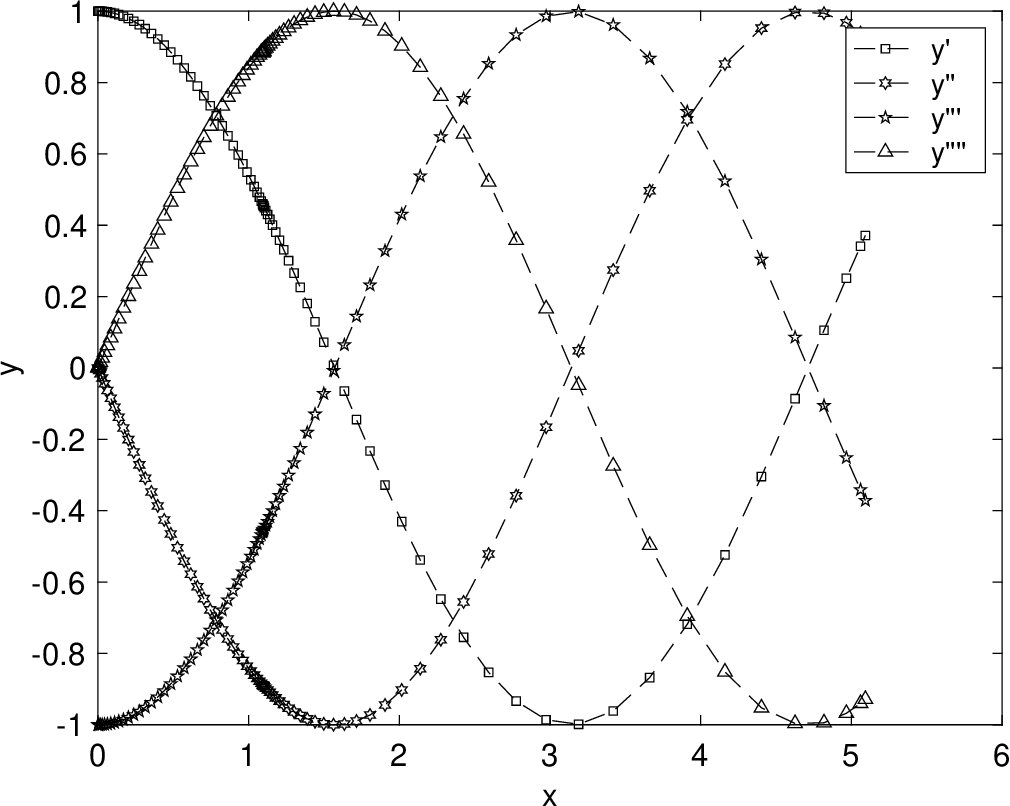}
\caption{High order derivative of $sin(x)$ for two and five points overlapped method: (a) two points  (b) five points.}
\label{y_yp4}
\end{figure}

In \Cref{y_yp4}, we show the difference between two and five-point overlapped techniques for fourth-order derivative and we see that the case with two points overlapped is overshot around the overlapped region. However, this bug has been fixed by increasing the points to five and more as shown in  \Cref{y_yp4}(b). 

\subsection{Algorithm 4: Multi points overlapping method}
While sections 2.2 and 2.3 showed how overlapping domains can be used to increase accuracy at the cost of slightly more points and a global solution matrix, it still requires constraints on the mappings, such as the need for coincident point(s) from the two subdomains.  Now, we do away with this requirement and consider arbitrary and independent grid point distributions in each subdomain, with only the requirement that the two domains overlap with at least two points from each subdomain in the overlap region. The basic premise is that points from one subdomain which lie within the interior of the other subdomain can be approximated using a polynomial expansion applied to the second domain   Thus, we use information from both $H_1$ and $H_2$ for each point located in the overlapping region. In principle, any form of quadrature or interpolation is possible; for example, one may consider using Lagrange interpolation to find the coefficients of points on $y_2$ in relation to $y_1$, or vice versa. However, this method introduces Runge's phenomenon despite the absence of equidistant points. Similarly, it is appealing to use Chebyshev interpolation to approximate the unknown function at points in the overlap region that belongs to one subdomain, by using the representation from the other domain; however, the matrix representation becomes quite complex. 

We propose combining Taylor's expansion with the Chebyshev spectral method to assemble the global matrix. Suppose we have a point $y_0$ located in the overlapping region $\delta$, which corresponds to Chebyshev points in layer $y_1$ but not in layer $y_2$. Now we need to estimate the derivative value at \( y_0 \) using the Chebyshev derivative matrix of layer \( y_2 \). We denote \( y_L \) as the closest point to the left of \( y_0 \), which is also a Chebyshev point in layer \( y_2 \), and $D_i$ is the $i-th$ order Chebyshev derivative matrix where $i = 1,2,3,...$. Then we estimate the derivative value at \( y_0 \) using information from layer $y_2$ as follows:

\begin{equation}
\begin{aligned}
u^{\prime}\left(y_0\right) 
& =u^{\prime}\left(y_L\right)+u^{\prime\prime}\left(y_L\right) \Delta y+\frac{u^{\prime \prime\prime}\left(y_L\right) \Delta y^2}{2!}+\frac{u^{\prime \prime \prime\prime}\left(y_L\right)}{3!} \Delta y^3 \\
& +\cdots \frac{u^{N+1}\left(y_L\right)}{N!}(\Delta y)^N+\frac{u^{N+2}(\xi)}{(N+1)!}(\Delta y)^{N+1}
\end{aligned}
\end{equation}

where $\Delta y = y_L-y_0$, and $u^{N}(y_L) = D_{N,t} * u_2(y)$, where $D_{N,t}$ is the $t$-th row of the $N$-th order Chebyshev derivative matrix. $\xi$ is some real number between $y_0$ and $y_L$. Then we obtain:

\begin{equation}
\begin{aligned}
u^{\prime}\left(y_0\right) 
& =u_2^{\prime}\left(y_L\right)+D_{2,t} u_2(y) \Delta y \,+\frac{D_{3,t} u_2(y) \Delta y^2}{2}+\frac{D_{4,t} u_2(y) \Delta y^3}{6} \\
& +\cdots \frac{D_{N,t} u_2(y) \Delta y^N}{N!}+\frac{D_{N+1,t}u_2(y)U(\xi) \Delta y^{N+1}}{(N+1)!} 
\end{aligned}
\end{equation}

Thus, we have obtained the derivative $u^{\prime}(y_0)$ for a point on $y_1$ by using the polynomial representation based on all Chebyshev nodes in layer $y_2$. As stated previously, the value of $u^{\prime}(y_0)$ can be found using the Chebyshev derivative matrix on $y_1$. Let $T_i$ denote the $i$-th order Chebyshev derivative matrix on layer $y_1$. Then one can get $u^{\prime}(y_0) = T_{1,s}u_1(y)$.  Here $s$ is the row number of $y_0$ in layer $y_1$ among its Chebyshev grid points. Next, if we assign weights of  $1/2$ to the equations resulting from these points (one from the Chebyshev representation and the other from Taylor's expansion based on the other layer) to estimate  $u^{\prime}(y_0)$, we  get 

\begin{equation}
\begin{aligned}
\left.\frac{d u}{d y}\right|_{y=y_0} 
& =  [u_2^{\prime}\left(y_L\right)+D_{2,t} u_2(y) \Delta y \,+\frac{D_{3,t} u_2(y) \Delta y^2}{2}+\frac{D_{4,t} u_2(y) \Delta y^3}{6}\\
& +\cdots \frac{D_{N,t} u_2(y) \Delta y^N}{N!} ] / 2 + T_{1,s}u_1(y)
\end{aligned}
\end{equation}

\begin{figure}
\centering
\includegraphics[height=2.0in]{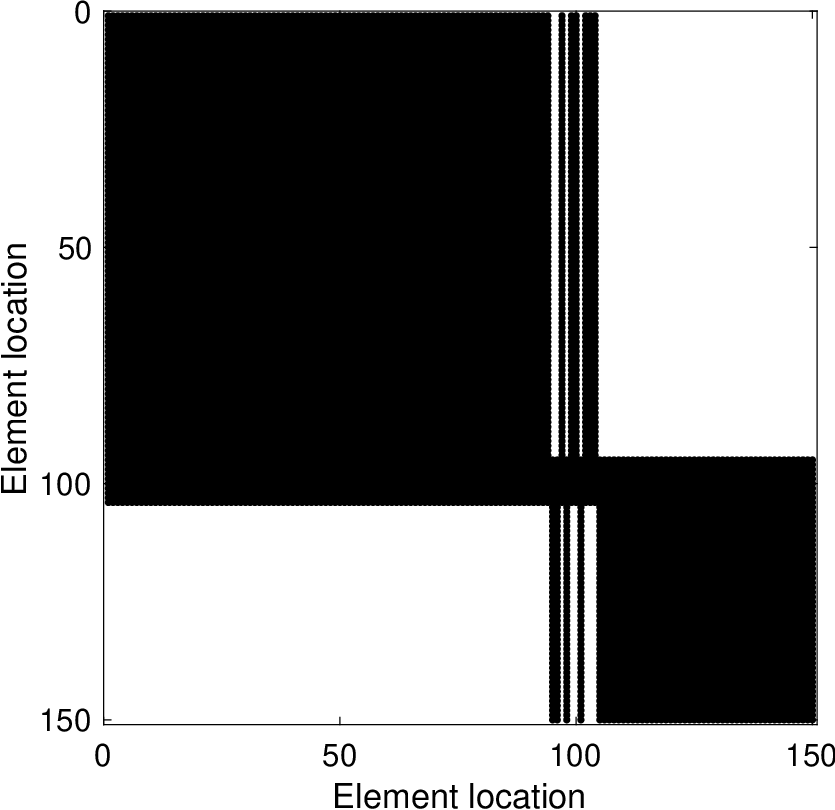}
\includegraphics[height=2.0in]{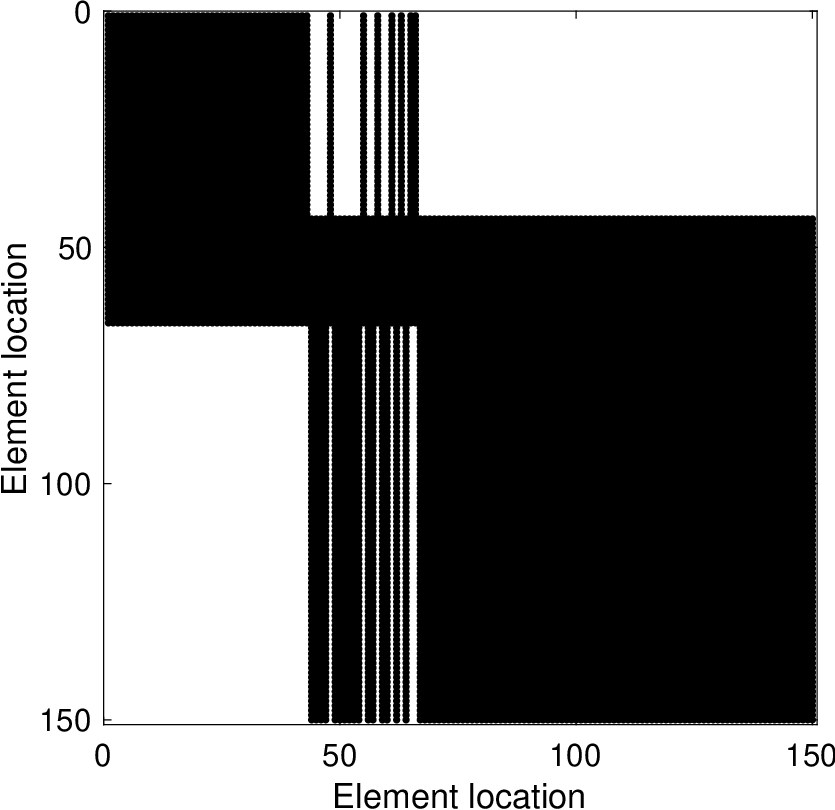}
\caption{Non-zero element distribution for (a) ten and (b) thirty points asymmetrical overlapping method.}
\label{fig:y_yprime2}
\end{figure}

\begin{figure}
\centering
\includegraphics[width=0.49\textwidth]{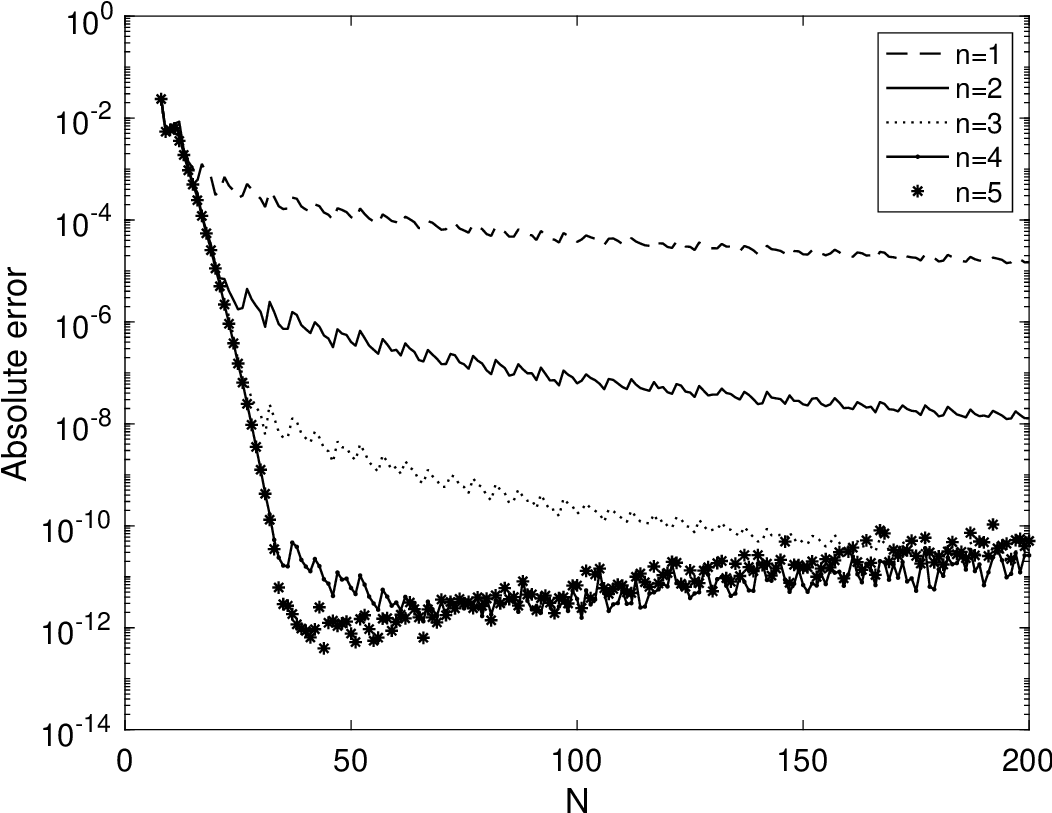}
\includegraphics[width=0.49\textwidth]{Figs/abs_err_mul_du.eps}
\caption{Convergence rate for multi-points overlapped method for Burgers equation with different terms of the Taylor expansion (n=1,2,3,4) for (a): two points overlapping and (b): five points overlapping.}
\label{fig:y_convg}
\end{figure}

Now we can see that \( u'(y_0) \) is estimated globally. For all points clustered in the overlapping region, we apply the same procedure, which allows us to obtain the global derivative matrix. Figure \ref{fig:y_yprime2} shows a schematic of the derivative and second derivative matrixes constructed using this algorithm, for 150 points in the two domains, with 10 points on the overlap region. Figure \ref{fig:y_convg} shows the convergence behavior for solutions of the Burgers equation, with the same parameter choices used previously. It may be observed from Fig. \ref{fig:y_convg}(a) that no more than 4 terms of the Taylor expansion are required to ensure exponential convergence of the solution. Similarly, Fig. \ref{fig:y_convg}(b) shows that no more than four points are required in the overlap region to ensure exponential convergence, which shows the number of overlapping points doesn't have any effect on the convergence rate but the terms of Taylor expansion make a big difference. 

\subsection{ Multi points overlapping on multi-interval:} With fewer grid points in each subdomain, a more computationally efficient method can be achieved. Over the past decade, the two-grid-point overlapping method has been extensively studied, requiring the use of linear mapping to ensure the two overlapped grid points are equidistant from both sides. However, our multi-point overlapping method does not have such constraints, allowing for the introduction of both linear and nonlinear mapping to expand the overlapping grid points. In the following section, we will provide the details of the numerical procedure for the multi-point overlapping method on multiple intervals.

Let's focus on the domain $I = [a, b]$ and we split it into $N$ subintervals and denote as $I_i, I = 1,2,3...N$. Here we can just consider using linear mapping to simplify the procedure but our method also works for nonlinear mapping. So each subinterval is $\Delta = (b-a)/N$  Now we consider the overlapping region to be $\delta$ for each subinterval. Then we have

\begin{equation}
\begin{gathered}
I_1=[a, a+\Delta+\delta] \\
I_i=[a+(i-1) \Delta-\delta, a+i \Delta+\delta]\\
I_N=[b-\Delta-\delta, b]
\end{gathered}
\end{equation}
where $i=2,3, \ldots, N-1$. We introduced the Chebyshev points,
$x_j=\cos (j \pi / N), \quad j=0,1, \ldots, N$. $D_x$ and $D_{xx}$ are the first and second-order derivative matrix for $x$ defined on $[-1,1]$. We can denote $I_{ai}$ and $I_{bi}$ as the beginning and end points of each subinterval. Then the Chebyshev differentiation matrices for each domain can be expressed as 
\begin{equation}
\begin{gathered}
D_i=\frac{2 D_x}{I_{a i}-I_{b i}} \\
D_i^j= (\frac{2 D_x}{I_{a i}-I_{b i}})^j
\end{gathered}
\end{equation}
and
\begin{equation}
x_i= \frac{I_{ai}-I_{bi}}{2}x +  \frac{I_{ai}+I_{bi}}{2}
\end{equation}

\begin{figure}
\centering
\includegraphics[height=2.0in]{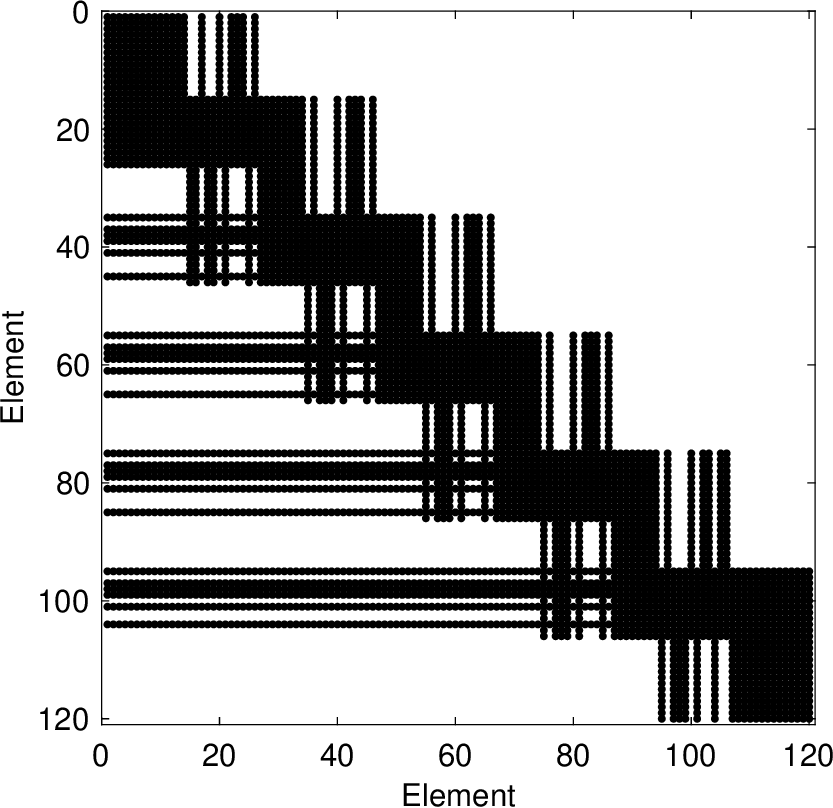}
\includegraphics[height=2.0in]{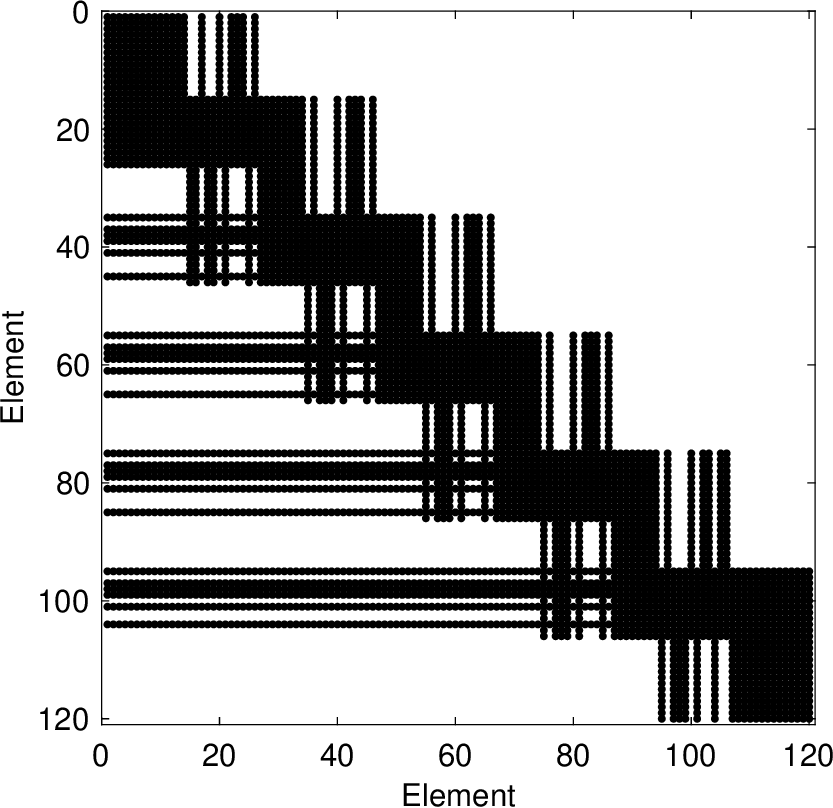}
\caption{The no-zero element distribution for the first and third order derivative matrix of the multi-point overlapped method with a six-layer interval and each layer with 20 nodes. The non-zero element distributions for any order of derivative are the same since they are assembled via a similarity algorithm.}
\label{yp_yp3}
\end{figure}


\clearpage
\section{Conclusions}
\label{sec:conclusions}
This paper presents four novel domain decomposition algorithms integrated with nonlinear mapping techniques to address collocation-based solutions of eigenvalue problems involving sharp interfaces or steep gradients. The proposed methods leverage the spectral accuracy of Chebyshev polynomials while overcoming limitations of existing tools like Chebfun, particularly in preserving higher-order derivative continuity and enabling flexible node clustering near discontinuities. Key findings include the following:  for algorithm Performance: The one-point overlap method demonstrated significant improvements over global mapping approaches, reducing required grid points by orders of magnitude (e.g., 300 vs. 6500 nodes) while maintaining spectral convergence. The two-point overlap method further enhanced accuracy by incorporating information from adjacent subdomains, resolving discontinuities in higher-order derivatives observed in Chebfun’s piecewise representations. The multi-point overlap methods generalized the approach, allowing arbitrary node distributions and nonlinear mappings. These achieved exponential error reduction (e.g., peak errors of 1e-8) for Burgers’ equation) by combining Taylor expansions with Chebyshev derivatives in overlap regions. While Chebfun’s splitting strategy automates domain decomposition, it enforces only $C_0$ continuity, leading to discontinuous higher derivatives. In contrast, the proposed algorithms preserved smoothness up to $C_N$ continuity, critical for eigenvalue problems in hydrodynamic stability and nonlinear BVPs. Validation on 3D channel flow with viscosity stratification and Burgers’ equation highlighted the methods’ robustness. For instance, eigenvalue calculations for miscible core-annular flows matched prior results while resolving sharp viscosity gradients with fewer nodes.

\newpage
\appendix
\section{Assembly global matrix of multi-point overlapping method1}

\begin{lstlisting}
shift_ind = 5; 
[x1,DM1] =chebdif(N0(1), 2); 
Dx1=DM1(:,:,1);
[x2,DM2] = chebdif(N0(2), 2);
Dx2=DM2(:,:,1);
y1 = (0-Rmd_up)/2*x1 + (Rmd_up+0)/2;
yy1 = y1;
dx1dy = 2/(0-Rmd_up);
Dy1 = dx1dy.*Dx1;
Dyy1 = dx1dy.*Dx1;
b = 1.25;
L = 2;
y2 = L*(1-x2)./(b+x2)+y1(N0(1)-shift_ind);
yy2 = y2;
Dy2 = ((b+x2).^2./(-L*(b+1))).*Dx2;
Dyy2 = ((b+x2).^2./(-L*(b+1))).*Dx2;
for jj = 1:length(y2)
for ii = 1:length(y1)-1
    if y2(jj) > y1(ii) AND y2(jj) <  y1(ii+1)
        ind_y1 = ii;
        Dy1 = [Dy1(:,1:ind_y1),zeros(N0(1),1),Dy1(:,ind_y1+1:end)];
        y1 = [y1(1:ind_y1);y2(jj);y1(ind_y1+1:end)];
    else
    end
end
end
Dy1 = [Dy1(:,1:N0(1)-shift_ind),zeros((N0(1)),1),...
Dy1(:,N0(1)-shift_ind+1:end)];
D = [Dyy1;Dyy2];
r = [yy1;yy2];
[a,b] = sort(r,'ascend'); % b is the reorder number based on the output a
for ii = 1:length(b)
    r1(ii,1) = a(ii);
    D1(ii,:) = D(b(ii),:);
end
DDy1(1,:) = [Dy1(1,:),zeros(1,N0(1)+N0(2)-length(Dy1(1,:)))]; 
%  the first element for global matrix; only this is enough
%  a template to arrange the points distribution
for tt = 1:N0(1)+N0(2)
    if ismember(r1(tt),yy1) ==1
        ind_add = 0;
        for jj = 1:N0(1)+N0(2)
            if logical(DDy1(1,jj)) == 1
                DD(tt,jj) = D1(tt,jj-ind_add);
            else
                DD(tt,jj) = 0;
                ind_add = ind_add+1;
            end
        end
    else
        ind_add = 0;
        for jj = 1:N0(1)+N0(2)
            if logical(DDy1(1,jj)) == 0
                DD(tt,jj) = D1(tt,jj-ind_add);
            else
                DD(tt,jj) = 0;
                ind_add = ind_add+1;
            end
        end
    end
end
ind_add = 0;
for jj = 1:N0(1)+N0(2)
    if logical(DDy1(1,jj)) == 0
        DD(N0(1)-shift_ind+1,jj) = D1(N0(1)-shift_ind+1,jj-ind_add);
    else
        DD(N0(1)-shift_ind+1,jj) = 0;
        ind_add = ind_add+1;
    end
end
D_st = [DD(:,1:N0(1)-shift_ind-1),DD(:,N0(1)-...
    shift_ind)+DD(:,N0(1)-shift_ind+1),DD(:,N0(1)-shift_ind+2:end)];
D_st(N0(1)-shift_ind,:) = D_st(N0(1)-shift_ind,:)./2+...
    D_st(N0(1)-shift_ind+1,:)./2;
D_st(N0(1)-shift_ind+1,:) = [];
r1(N0(1)-shift_ind+1) = [];
d1 = D_st;
d2 = D_st*D_st;
Dr = D_st;
Drr = D_st*D_st;
\end{lstlisting}

\section{Assembly global matrix of multi-point overlapping method 2.}

\begin{lstlisting}
N0 = [60,30];
theta_mu = 0.1;
loc = 1;
Rmd_up = loc+theta_mu/2;
[x1,DM1] =chebdif(N0(1), 2);
Dx1=DM1(:,:,1);
[x2,DM2] = chebdif(N0(2), 2);
Dx2=DM2(:,:,1);
y1 = (0-Rmd_up)/2*x1 + (Rmd_up+0)/2;
yy1 = y1;
dx1dy = 2/(0-Rmd_up);
Dy1 = dx1dy.*Dx1;
Dyy1 = Dy1*Dy1;
Dyyy1 = Dy1*Dyy1;
Dyyyy1 = Dyy1*Dyy1;
Dyyyyy1 = Dyyy1*Dyy1;
Dyyyyyy1 = Dyyy1*Dyyy1;
b = 3;
L =  1.5;
y2 = L*(1-x2)./(b+x2)+y1(end)-0.1/2;
yy2 = y2;
Dy2 = ((b+x2).^2./(-L*(b+1))).*Dx2;
Dyy2 = Dy2*Dy2;
Dyyy2 = Dyy2*Dy2;
Dyyyy2 = Dyyy2*Dy2;
Dyyyyy2 = Dyyyy2*Dy2;
Dyyyyyy2 = Dyyyyy2*Dy2;
ri = [yy1;yy2];
ri = sort(ri);
indices_y1 = ismember(ri, y1);
log_y1 = find(indices_y1);
indices_y2 = ismember(ri, y2);
log_y2 = find(indices_y2);
Dr = zeros(N0(1)+N0(2));
for ii = 1:length(ri)
    if ri(ii) >= y1(1) && ri(ii) < y2(1) 
       Dr(ii,log_y1)  = Dy1(ii,:);
    elseif ismember(ri(ii),y1)  && ismember(ri(ii),y2) 
        index_a = find(y1 == ri(ii));
        index_b = find(y2 == ri(ii));
       Dr(ii,log_y1)  = Dy1(index_a,:)./2;
       Dr(ii,log_y2) = Dy2(index_b,:)./2;
    elseif ri(ii) >= y2(1) && ri(ii) <= y1(end) &&  ismember(ri(ii),y1) 
     index_1 = find(y1 == ri(ii));
     index_left = max((find(ri(ii) > y2)));  % find the most right on y2 for ri(ii)
      Dr(ii,log_y1) = Dy1(index_1,:)./2;
      delta_x = ri(ii)-y2(index_left);
      Dr(ii,log_y2) = ( Dy2(index_left,:) + delta_x*Dyy2(index_left,:) + delta_x^2./2*Dyyy2(index_left,:) ...
          + delta_x^3./6*Dyyyy2(index_left,:) ...
          + delta_x^4./24*Dyyyyy2(index_left,:) ...
          + delta_x^5./24/5*Dyyyyyy2(index_left,:) )./2;  
    elseif ri(ii) >= y2(1) && ri(ii) <= y1(end) &&  ismember(ri(ii),y2) 
     index_2 = find(y2 == ri(ii));
     index_left_2 = max((find(ri(ii) > y1))); 
     delta_x2 = ri(ii)-y1(index_left_2);   % ri(ii) on y2 so we need to use Dy1 
       Dr(ii,log_y1) = ( Dy1(index_left_2,:) + delta_x2*Dyy1(index_left_2,:) + delta_x2^2./2*Dyyy1(index_left_2,:) ...
          + delta_x2^3./6*Dyyyy1(index_left_2,:) + ...
          + delta_x2^4./24*Dyyyyy1(index_left_2,:) ...
          + delta_x2^5./24/5*Dyyyyyy1(index_left_2,:))./2;
      Dr(ii,log_y2)= Dy2(index_2,:)./2;
    elseif ri(ii) > y1(end)
        add5 = add5 + 1;
         index_3 = find(y2 == ri(ii));
          Dr(ii,log_y2) = Dy2(index_3,:);
    else
        add6 = add6+1;
    end
end
Drr = Dr*Dr;
z = ri;
\end{lstlisting}

\section*{Acknowledgments}
We want to acknowledge the NSF funding.

\bibliographystyle{siamplain}
\bibliography{references}
\clearpage

\end{document}